\documentclass[journal,12pt,twocolumn]{IEEEtran}

\usepackage{soul}
\usepackage{comment}
\usepackage{amssymb}
\usepackage{srcltx}
\usepackage{amsmath}
\usepackage{amsfonts}
\usepackage{psfrag}
\usepackage{graphicx,subfigure}
\usepackage{multirow}
\usepackage{algorithm}
\usepackage{algorithmic}   
\usepackage{psfrag}       
\usepackage{amsmath}    
\usepackage{amsfonts}
\usepackage{url}
\usepackage{color}
\usepackage{xcolor}
\usepackage{ulem}
\usepackage{caption}
\usepackage{subfigmat}

\graphicspath{ {figures/} }

\begin{document}

\title{Data-Driven Probabilistic Methodology for Aircraft Conflict Detection Under Wind Uncertainty}

\author{Jaime de la Mota, Mar\'{\i}a Cerezo-Maga\~na, Alberto Olivares, Ernesto Staffetti
\thanks{{Universidad Rey Juan Carlos, Department of Telecommunication Engineering, Fuenlabrada, Madrid, Spain (e-mail: ernesto.staffetti@urjc.es (Corresponding Author: Ernesto Staffetti).  This work was supported in part by the Spanish Government under Grant RTI2018-098471-B-C33 and Grant PID2021-122323OB-C31.}}}

\markboth{DE LA MOTA ET AL.}{DATA-DRIVEN PROBABILISTIC METHODOLOGY FOR AIRCRAFT CONFLICT DETECTION}
\maketitle

\begin{abstract}
Assuming the availability of a reliable aircraft trajectory planner, this paper presents a probabilistic methodology to detect conflicts between aircraft in the cruise phase of flight in the presence of wind velocity forecasting uncertainty. This uncertainty is quantified by ensemble weather forecasts, the members of which are regarded as realizations of correlated random processes and used to derive the eastward and northward components of the wind velocity. First, the Karhunen-Lo\`eve expansion is used to obtain a series expansion of the components of the wind velocity in terms of a set of uncorrelated random variables and deterministic coefficients. Then, the uncertainty generated by these uncorrelated random variables in the outputs of the aircraft trajectory planner is quantified using the arbitrary polynomial chaos technique. Finally, the probability density function of the great circle distance between each pair of aircraft is derived from the polynomial expansions using a Gaussian kernel density estimator and used to estimate the probability of conflict. The arbitrary polynomial chaos technique allows the effects of uncertainty in complex nonlinear dynamical systems, such as those underlying aircraft trajectory planners, to be quantified with high computational efficiency, only requiring the existence of a finite number of statistical moments of the random variables of the Karhunen-Lo\`eve expansion while avoiding any assumptions on their probability distributions. To demonstrate the effectiveness of the proposed conflict detection method, numerical experiments are conducted via an optimal control-based aircraft trajectory planner for a given wind velocity forecast represented by an ensemble prediction system.
\end{abstract}


\section{INTRODUCTION}
\label{sect:intro}


During recent decades, air traffic has grown at a steady rate. This growth has been hampered by several crises, though, historic data show that commercial aviation has always recovered. To cope with the increase in air traffic, in 2005, the International Civil Aviation Organization (ICAO) proposed a new operational concept for Air Traffic Management (ATM) \cite{andreeva2013rule}. In the European Union, the Single European Sky ATM Research (SESAR) \cite{undertaking2015european} is an ongoing program aimed at developing new ATM systems to improve the coordination of air traffic operations, thus increasing the efficiency, capacity, and safety of the current ATM system, while decreasing the environmental impact of commercial aviation. Similar programs are being developed in other countries, such as the Next Generation Air Transportation System (NextGen) \cite{planning2010concept}, in the USA, and the Collaborative Actions for Renovation of Air Traffic Systems (CARATS), in Japan \cite{bureau2010long}.
The main objective of the ATM is to ensure the flight safety. 


Aircraft conflict is one of the most dangerous situations in air traffic. It occurs when the distance between two or more aircraft is below the required minimum separation. More specifically, according to the current regulations \cite{icao4444doc}, the horizontal separation distance should be at least \mbox{$5$ [NM]}, whereas the vertical separation distance should be at least 1000 [ft].

Aircraft conflict avoidance systems operate in two stages: detection and resolution. In the detection stage, aircraft trajectories are predicted in the relevant airspace. Then, conflict is identified when the probability of losing separation distance exceeds a certain threshold. In the resolution stage, after selecting the most suitable maneuvers, the predicted trajectories are adjusted to prevent the detected conflict. The sooner a conflict is detected, the more efficient the Conflict Resolution (CR) maneuvers in the resolution phase.
Moreover, early Conflict Detection (CD) also reduces the workload of the air traffic controllers.


There are various sources of uncertainty that influence the accuracy of the aircraft trajectory prediction, such as data uncertainty, operational uncertainty (the lack of knowledge in the decisions taken by individuals), equipment uncertainty (associated with malfunctions and breakdowns of communication, navigation, and surveillance systems), and weather uncertainty, which occurs when meteorological phenomena, mainly wind and thunderstorms, are not accurately forecast. Weather uncertainty in particular affects air traffic as a whole and has a very strong influence on the accuracy of the trajectory prediction \cite{cook2016complexity}.


This paper investigates the aircraft 
CD problem in the presence of wind velocity forecasting uncertainty. More specifically, the CD problem is studied in the cruise phase of flight, in which probabilistic wind velocity forecasts generated by Ensemble Prediction Systems (EPS) are used to represent wind velocity uncertainty.

Weather forecasts are estimates of the future state of the atmosphere. They are usually generated by estimating the initial conditions of the atmosphere using observations and then computing the temporal evolution using a Numerical Weather Prediction (NWP) model. However, since these NWP models are very sensitive to the initial conditions, small errors in the initial state may result in large errors in the forecast. 
In EPS forecasts, the NWP model is run several times using slightly different initial conditions and model parameters, and then a set of forecasts is obtained, which is referred to as the ensemble. The individual forecasts of the ensemble are called members, which are assumed to be equiprobable. The differences between their corresponding initial states are consistent with the uncertainty in the observations. The goal of EPS, which typically are composed of between 10 and 50 members, is to represent the various possible evolutions of the atmospheric conditions. Ideally, the future state of the atmosphere should be within the limit of the spread of the predicted ensemble. Further details on EPS can be found in \cite{palmer2019ecmwf, buizza2005comparison, molteni1996ecmwf}.


The essence of the probabilistic approach to the aircraft CD problem is to estimate the probability of loss of separation between aircraft. The aircraft CD problem has been traditionally associated with the aircraft 
CR problem. Aircraft CD and Resolution (CD\&R) methods can be classified, according to the look-ahead time at which they operate, into strategical, tactical, and operational, which operate with look-ahead times of more than 30 [min], between 30 and 10 [min], and less than 10 [min], respectively.
There is a recent survey paper on aircraft CD\&R,  \cite{tang2019conflict}, which includes literature on aircraft CD in the presence of wind velocity forecasting uncertainty. Another recent survey paper on stochastic modeling in aviation, \cite{shone2021applications}, also includes literature on aircraft CD in the presence of wind velocity forecasting uncertainty.

In aircraft CD, the aircraft trajectories are first predicted and then analyzed to calculate the probability of conflict between them. In particular, in aircraft CD in which the wind velocity forecast is represented by EPS, the uncertainty contained in the ensemble forecast has to be propagated to the prediction of the trajectories. There are two common approaches to trajectory prediction for CD using wind velocity forecasts generated by EPS: the transformation approach and the ensemble approach.

 In the transformation approach, the probability distributions of the relevant uncertain meteorological parameters, such as the eastward and northward components of the wind velocity, are estimated from the ensemble forecast. Then, a probabilistic aircraft trajectory planner is used to obtain the probability distributions of the variables used to detect conflicts, such as the aircraft positions with respect to time.  In the ensemble approach, a deterministic aircraft trajectory planner is used to calculate a trajectory for each member of the ensemble wind velocity forecast, thus obtaining an ensemble of aircraft trajectories, from which  the probability distributions of the variables used to detect conflicts are derived.
Since it is less demanding from a computational point of view and provides similar results to the transformation approach \cite{rivas2017analysis}, the ensemble approach is more suitable for practical applications.

 \cite{hernandez2019probabilistic} presents a methodology to quantify statistically the severity of conflicts between aircraft flying at the same altitude considering the wind velocity forecasting uncertainty derived from the PEARP EPS released by M\'et\'eo-France, consisting of 35 members. The severity of the conflict is characterized by two descriptors: conflict intensity and conflict probability. Conflict intensity is defined as the mean minimum distance between a pair of aircraft approaching a waypoint. The components of the wind velocity are modeled as four-parameter beta distributions. The probability of conflict is obtained in terms of the Probability Density Function (PDF) of the minimum distance between aircraft, calculated from the PDFs of the components of the wind velocity using the transformation approach, which allows the joint PDF of a set of random variables to be obtained from the joint PDF of another set of random variables. \cite{hernandez2020probabilistic} combines the aircraft CD methodology proposed in  \cite{hernandez2019probabilistic} with a CR technique in which conflicts between aircraft are solved while minimizing the deviations from the original aircraft trajectories.

\cite{eulaliaTesis} presents a methodology for calculating the probability of conflict between aircraft flying 3D routes considering wind velocity and temperature forecasting uncertainty derived from the COSMO-D2 EPS released by the Deutscher Wetterdienst, consisting of 20 members. 
The ensemble approach is used, in which aircraft trajectories are calculated deterministically for each member of the ensemble \cite{gonzalez2018robust}, and then a deterministic CD is conducted for all the computed trajectories. Since all the members are assumed to be equiprobable, the probability of conflict between each pair of aircraft is calculated as the fraction of members for which a conflict between these two aircraft is identified.

\cite{jilkov2018multiple} proposes a multiple model method for aircraft  CD\&R considering intent and weather uncertainty.  It is based on probabilistic multiple model aircraft trajectory prediction. 
If a multiple-model trajectory prediction is used, the separation vector between two aircraft has a Gaussian mixture distribution, and an efficient randomized algorithm is proposed to estimate the probability of conflict.

\cite{Mishra2019} presents an efficient method for estimating the probability of conflict between aircraft, which is based on what is known as the subset simulation technique, in which the small conflict probabilities are computed as the product of larger conditional conflict probabilities, reducing the computational cost and improving the accuracy of the estimated probability.

This paper presents a methodology to compute the probability of conflict between aircraft flying at the same altitude considering wind velocity forecasting uncertainty quantified by ensemble weather forecasts. More specifically, EPS by the European Center for Medium-range Weather Forecast (ECMWF), consisting of 51 members, have been considered.

The eastward and northward components of the wind velocity have been modelled as two correlated random processes and the different members of the ensemble have been regarded as realisations of the two random processes. The multiple uncorrelated Karhunen-Lo\`eve (muKL) expansion \cite{cho2013karhunen} has been applied to the ensemble to reduce the dimension of the representation of the two correlated random processes. The muKL expansion extracts a set of uncorrelated random variables, which represent the uncertainty contained in the random processes. These random variables then form the random inputs of the optimal control model used to generate the aircraft trajectories. The resulting model is a complex stochastic optimal control model, which is difficult to solve. Therefore, a surrogate model has been devised, which is computationally easier to solve and accurately approximates the propagation of the uncertainty - represented by the uncorrelated random variables - on the state variables of the solution of the optimal control problem. The state variables of this solution include the optimal trajectories of the aircraft.
Polynomial Chaos Expansion (PCE) \cite{xiu2010numerical} has been used to formulate this surrogate model. PCE is a spectral method, which consists in the projection of the model outputs on a basis of orthogonal stochastic polynomials in the random inputs. This probabilistic method yields an efficient representation of the variability of the model outputs with respect to the variability of its inputs. In particular, in this paper, a data-driven moment-based PCE technique called arbitrary Polynomial Chaos (aPC) has been employed \cite{ahlfeld2016samba}.


The Karhunen-Lo\`eve (KL) expansion is one of the most common approaches for dimension reduction in the representation of random processes as the infinite sum of orthogonal deterministic basis functions, called eigenfunctions, multiplied by uncorrelated random coefficients. This infinite sum can be truncated depending on the precision required in the representation of the random process. The KL expansion has been developed for a single random process or ensembles of statistically independent random processes, and, therefore, its generalization to multi-correlated processes, namely to several random processes with mutual correlation, is not straightforward. Indeed, if the cross-covariances of several random processes are not zero, it is not easy to calculate consistent expansions for all the random processes that reflect both the structure of the autocorrelation and the structure of the cross-covariance. To overcome this difficulty, the muKL expansion has been introduced in \cite{cho2013karhunen}, which extends the classical KL expansion to multi-correlated nonstationary random processes and is based on the spectral decomposition of a suitable assembled random process, providing series expansions of the random processes using a single set of uncorrelated random variables.


The aPC expansion is a statistical moment-based PCE technique that allows surrogate models to be built. It extends the PCE techniques that require knowledge of the PDF of the input random variables, such as the generalized polynomial chaos expansion  \cite{xiu2010numerical}. Indeed, the aPC expansion only requires the existence of a finite number of moments of the input random variables of the model, which can be discrete or continuous and can be specified analytically as PDFs or numerically as histograms or raw data sets. This means that this polynomial expansion does not require the existence of a parametric PDF, avoiding the need to fit parametric probability distributions to data, which is very useful, especially when a limited amount of data is available.


In this paper, using the muKL expansion, the random processes that represent the components of the wind velocity are approximated by a finite sum of the product of certain eigenfunctions multiplied by their corresponding uncorrelated random variables, which are characterized using the aPC expansion.
%
%
%
%
First, the uncorrelated random variables obtained in the muKL expansion are ordered according to the amount of variability explained by each of them, which enables the infinite sum to be truncated according to the required precision. 
Then, the aPC expansion is calculated using the most significant random variables. A set of nodes and weights are computed in order to develop the surrogate models for the state variables of the optimal control model used to generate the aircraft trajectories, which represent the variability of the optimal trajectories of each aircraft.  
Next, the main statistics of the orthodromic distance between pairs of aircraft are estimated, from which the confidence envelopes are determined. 
Finally, the marginal and joint PDFs of the distance between aircraft at a given instant in time are estimated, and the marginal and conditional probabilities of conflict are computed. 
One of the main advantages of the data-driven methodology proposed in this paper is that it does not rely on any assumptions about these PDFs.


Although in this paper this technique has been applied to CD among aircraft, 
it can be seamlessly applied to CD among unmanned aerial vehicles, which are more sensitive than aircraft to wind velocity uncertainty and for which 
conflict resolution may entail solving complex cooperative multiple task reallocation problems 
\cite{tangetal:2022:drmomuavtieas}.


This paper is organized as follows:
Section~\ref{section:muKL}  presents the general formulation of the
muKL expansion, which is used to reduce the dimension of the representation of multi-correlated random processes.
Section~\ref{section:aPC} describes the moment-based PCE technique, which is employed for the uncertainty quantification.
Section~\ref{section:UQ} outlines the computational and statistical properties of the surrogate models  represented by the PCE expansion, which are used to estimate the marginal and joint PDFs of the distance between aircraft.
The practical application of the proposed methodology for aircraft CD is shown in Section~\ref{section:numerical_results}.
Finally, some conclusions are drawn in Section~\ref{section:conclusions}.

\section{Karhunen-Lo\`eve Expansion for Multi-correlated Random Processes}
 \label{section:muKL}


As mentioned before, the eastward and northward components of the wind velocity forecast represented by the EPS  are considered realizations of correlated random processes. Therefore, to implement their KL expansion, a specific methodology for correlated processes is used. In particular, the muKL expansion is used, which is one of the two methodologies proposed in \cite{cho2013karhunen} with the aim to generalize the KL expansion in order to model multi-correlated non-stationary random processes. The other approach is referred to as the Multiple Correlated KL (mcKL) expansion. The muKL technique has been chosen in this paper since it allows the series expansions of the correlated processes to be generated in terms  of a single set of uncorrelated random variables, whereas in the mcKL approach, a different set of mutually correlated random variables must be developed for each process. In this section, the muKL method is described for convenience. A more in-depth treatment of the method, accompanied by illustrative examples, is given in \cite{cho2013karhunen}.

Let $(\Omega, \Lambda, P)$ be a probability space, where $\Omega$  is the sample space, $\Lambda$ is a $\sigma$-algebra, and P is a probability measure. Consider the following ensemble of $n$ zero-mean, square-integrable random processes:

\begin{equation}
\{f_1(s;\omega), \ldots, f_n(s; \omega))\}, \;\;\; \omega \in \Omega.
\label{eq:karhunen-loeve_1}
\end{equation}
Without loss of generality, it is assumed that all these processes are defined in the same bounded interval $[0, S]$, where index $s$ usually denotes time or space.
More specifically, in this paper, index $s$ represents the spatial position in which weather data are predicted. Then, the correlation between the aforementioned processes is defined in terms of the following $n(n + 1)/2$ covariance kernels:

\begin{multline}
C_{ij}(s_1,s_2) = \\
E[f_i(s_1; \omega) f_j(s_2; \omega)], \; 1 \; \leq i \leq j \leq \; n,
\label{eq:covariance_functions}
\end{multline}
where $E[ \; \cdot \; ]$ is the statistical expectation operator. In particular, $C_i(s_1, s_2)=C_{ii}(s_1, s_2)$ represents the auto-covariance function of the process $f_i(s, \omega)$.

Notice that if the $n$ processes (\ref{eq:karhunen-loeve_1}) are mutually independent, then the conventional KL expansion can be directly applied to each process, resulting in multiple series which can be built separately \cite{huang2001convergence}. Conversely, if the cross-covariances (\ref{eq:covariance_functions}) are not zero, then the classical KL expansion is not capable of providing consistent expansions for all the random processes, which reflect the autocorrelation and the cross-covariance structure.

The muKL approach overcomes this drawback by generating a series expansion for each random process of the ensemble processes (\ref{eq:karhunen-loeve_1}) in terms of a single set of uncorrelated random variables \cite{ramsay2002applied}. To obtain such a series, the following assembled process is defined:
\begin{equation}
\tilde{f}(s, \omega) = f_i(s-S_{i-1}; \omega), \;\;\;\; s \in \mathcal{I}_i,
\label{eq:paper_cho_2}
\end{equation}
where $S_i=iS$ and $\mathcal{I}_i=(S_{i-1}, S_i], 1\leq i \leq n$, with $\mathcal{I}_1=[0, S_1]$. That is to say, the restriction of the assembled process  $\tilde{f}(s, \omega)$ to the interval $\mathcal{I}_i$ corresponds to the process $f_i(s, \omega)$. Notice that  $\tilde{f}(s, \omega)$ is a second-order process as well, which satisfies

\begin{equation*}
E[\tilde{f}(s, \omega)]=0, \;\; E[\tilde{f}(s_1, \omega)\tilde{f}(s_2, \omega)]= \tilde{C}(s_1,s_2),
\end{equation*}   
where $\tilde{C}(s_1, s_2)$ is the assemble covariance function, which is defined as
\begin{multline}
\tilde{C}(s_1, s_2) = \\ C_{ij}(s_1-S_{i-1}, s_2-S_{j-1}),   \;  s_1 \in \mathcal{I}_i,  \;  s_2 \in \mathcal{I}_j.
\label{eq:paper_cho_3}
\end{multline}

Then, a conventional KL expansion can be applied to the assembled process (\ref{eq:paper_cho_2}) obtaining
\begin{equation}
\tilde{f}(s, \omega)= \sum_{k=1}^{\infty} \sqrt{\lambda_k}\tilde{f}_k(s)\xi_k(\omega),
\label{eq:paper_cho_4}
\end{equation}
with $\xi_k(\omega)$ uncorrelated random variables, which are calculated as
\begin{equation*}
\xi_k(\omega) = \frac{1}{\sqrt{\lambda_k}}\int_{0}^{S_n} \tilde{f}(s, \omega) \tilde{f}_k(s) ds,
\end{equation*}
with $\lambda_k$ and $\tilde{f}_k(s)$ defined as  eigenvalues and eigenfunctions of a symmetric compact integral operator \cite{kato2013perturbation} whose kernel is (\ref{eq:paper_cho_3}), namely $\lambda_k$ and $\tilde{f}_k(s)$ are solutions to the following homogeneous Fredholm integral equation of the second kind:
\begin{equation}
\lambda_k \tilde{f}_k(s_1) = \int_{0}^{S_n} \tilde{C}(s_1, s_2)\tilde{f}_k(s_2) ds_2.
\label{eq:fredholm_2}
\end{equation}

In practical applications, it is more convenient to work with non-negative covariance functions. Unfortunately, the assembled covariance $\tilde{C}(s_1, s_2)$ might have negative eigenvalues due to the fact that, in general, it is not positive semi-definite, even if all the covariances (\ref{eq:covariance_functions}) are. Thus, a positivity condition for the assembled covariance $\tilde{C}(s_i, s_j)$ is imposed, that is to say
\begin{equation*}
 \displaystyle \sum_{j=1}^{m}  \sum_{i=1}^{m} \tilde{C}(s_i,s_j)x_ix_j \geq 0,
\end{equation*}
for any finite sequence $\{s_1, \ldots, s_m\}$ and any real numbers $x_i$, $1 \leq i \leq m$. Namely, the $m\times m$ matrix
 \[
   \tilde{C}=
  \left[ {\begin{array}{cccc}
   \tilde{C}(s_1, s_1) & \tilde{C}(s_1, s_2) &\cdots &\tilde{C}(s_1, s_m)  \\
   \tilde{C}(s_2, s_1) & \tilde{C}(s_2, s_2) &\cdots &\tilde{C}(s_2, s_m)  \\
   \vdots & \vdots &\ddots &\vdots  \\  
   \tilde{C}(s_m, s_1) & \tilde{C}(s_m, s_2) &\cdots &\tilde{C}(s_m, s_m)  \\
\end{array} } \right]
\]
must be positive semi-definite for any collection of $m$ different values of $s \in [0, S]$.

Thus, the eigen-pairs $\{ \lambda_k,  \tilde{f}_k(s)\}, k \geq 1$, can be calculated and arranged according to the magnitudes of the eigenvalues $\lambda_k$ using
(\ref{eq:fredholm_2}). Moreover, each eigenfunction $\tilde{f}_k(s)$ can be expressed in terms of $n$ subfunctions $\phi_k^{(i)}(s), 1 \leq i \leq n,$ which are defined as
\begin{equation*}
\phi_k^{(i)}(s) = \tilde{f}_k (s+S_{i-1})\mathcal{I}_{[0, S]}(s),
\end{equation*}
being $\mathcal{I}_{[0, S]}$ the indicator function in the interval $[0, S]$. Therefore, the $i$-th random process $f_i(s,\omega)$ of the original ensemble (\ref{eq:karhunen-loeve_1}) is represented as follows: 
\begin{equation}
f_i (s,\omega) =  \displaystyle \sum_{k=1}^{\infty} \sqrt{\lambda_k} \phi_k^{(i)}(s) \xi_k(\omega).
\label{eq:paper_cho_5}
\end{equation}

Notice that the eigenvalues $\lambda_k$ and the random variables $\xi_k(\omega)$ which appear in  (\ref{eq:paper_cho_5}), are the same as the ones introduced in (\ref{eq:paper_cho_4}). Moreover, for each index $i$, the collection of  subfunctions $\{ \phi_k^{(i)}(s)\}$, $ k \geq 1$,  is neither orthogonal nor normalized in $s \in[0, S]$. Nevertheless, this drawback can be overcome by normalizing $\phi_k^{(i)}(s)$ within $[0, S]$, namely the random process $f_i(s, \omega)$ is rewritten as
\begin{equation}
f_i(s, \omega) = \sum_{k=1}^{\infty} \sqrt{\hat{\lambda}_k^{(i)}} \hat{\phi}_k^{(i)}(s)\xi_k(\omega),
\label{eq:paper_cho_6}
\end{equation}
being $\hat{\phi}_k^{(i)}(s) = \phi_k^{(i)}(s)/\big\Vert\phi_k^{(i)}(s)\big\Vert_2$ and $\hat{\lambda}_k^{(s)}= \lambda_k\big\Vert\phi_k^{(i)}(s)\big\Vert_2^2$.

For practical purposes, the dimension of the expansion (\ref{eq:paper_cho_4}) is reduced via truncation, and its associated mean-square error can be calculated. A truncated assembled process, which only considers the first $M$ elements of the expansion, is defined as follows: 
\begin{equation}
Q_M(s, \omega) = \sum_{k=1}^{M}\sqrt{\lambda_k}\tilde{f}_k(s)\xi_k(\omega),
\label{eq:paper_cho_7}
\end{equation}
whereas the associated mean-squared error is calculated as
\begin{equation}
\varepsilon_M^2 = \int_0^{S_n} E[(\tilde{f}(s, \omega)- Q_M (s, \omega))^2]ds.
\label{eq:MeanSquareError}
\end{equation} 

Taking into account that $\xi_k(\omega)$ are uncorrelated random variables and $\tilde{f}_k(s)$ are orthonormal eigenfunctions, the mean-squared error (\ref{eq:MeanSquareError}) can be expressed as
\begin{equation}
\varepsilon_M^2= \sum_{k=M+1}^{\infty} \lambda_k.
\label{eq:paper_cho_8}
\end{equation} 
Thus, the truncation error of the series expansion (\ref{eq:paper_cho_4}) decreases with respect to the decay rate of the eigenvalues. Furthermore, the mean-squared error (\ref{eq:MeanSquareError}) is an upper bound on the truncation error of the series expansion (\ref{eq:paper_cho_6}), whereas the errors of the cross-covariances $C_{ij}(s_1, s_2)$ are bounded by the error of the assembled covariance  $\tilde{C}(s_1,s_2)$ defined in (\ref{eq:paper_cho_3}).

In practical applications, a key aspect of any KL expansion is the election of the appropriate number of terms in the truncated expansion. In particular, in the muKL expansion, a threshold for the error $\varepsilon_M$ can be imposed to determine the choice of $M$ in (\ref{eq:paper_cho_7}). According to (\ref{eq:paper_cho_8}), this can be achieved by imposing a threshold on the relative sum of eigenvalues, as follows:
\begin{equation*}
\sum_{k=1}^{M} \lambda_k \geq \delta \sum_{k=1}^{\infty} \lambda_k,
\end{equation*}
being $\delta \in [0, 1]$ an arbitrary constant chosen so that the accuracy of the approximation is satisfactory for a particular application. 

Notice that, in order to incorporate the muKL expansion into the aircraft trajectory planner model, the eigenfunctions $\phi_k^{(i)}(s), 1 \leq i \leq n$, associated with the expansions of the eastward and northward components of the wind velocity must be interpolated to be converted into analytic expressions.
Moreover, the uncertainty of the corresponding random variables $\xi_k(\omega)$ must be quantified, which can be done using the aPC expansion introduced in Section~\ref{section:aPC}.


\section{Moment-Based Arbitrary Polynomial Chaos}
 \label{section:aPC}

Following  \cite{ahlfeld2016samba}, this section introduces the aPC approach, which is used in this article to represent the propagation of the uncertainty through   the aircraft trajectory planner model.

Let $\boldsymbol{\xi}=(\xi_1, \xi_2, \ldots, \xi_{N_U})$ be a vector of $N_U$ independent random variables in the probability space $(\Omega, \Lambda, P)$ introduced in Section \ref{section:muKL}. Notice that, for the sake of simplicity, the formal dependency on $\omega$ is dropped for $\xi_i(\omega), 1 \leq i \leq N_U$. Then, a surrogate model for each output variable $x(t, \boldsymbol{\xi})$ of the aircraft  trajectory
planner can be computed, which is represented by a multidimensional polynomial expansion. In particular, a linear combination of $N_P$ stochastic multivariate orthonormal polynomials $\boldsymbol{\Psi}_k(\boldsymbol{\xi})$ with deterministic coefficients $\alpha_k(t)$ can be used to approximate the output variable $x(t, \boldsymbol{\xi})$ as follows:
\begin{multline}
x(t, \boldsymbol{\xi})=x(t; \xi_1,  \xi_2,\ldots, \xi_{N_U}) \approx \\
\sum_{k=1}^{N_P}\alpha_k(t)\boldsymbol{\Psi}_k(\xi_1, \xi_2, \ldots, \xi_{N_U}).
\label{eq:paper_alberto_3}
\end{multline}

In   (\ref{eq:paper_alberto_3}), the  multivariate orthonormal polynomials $\boldsymbol{\Psi}_k(\boldsymbol{\xi})$, with $1 \leq k \leq N_P$, are calculated as the product of univariate orthonormal polynomials $\psi_j^i(\xi_i)$,  with $1 \leq i \leq N_U$, $1 \leq j \leq p$, where $i$ indicates the elements of the vector of random variables  and $j$ the order of the univariate orthonormal polynomials. 
More specifically, each  polynomial $\boldsymbol{\Psi}_{k}(\boldsymbol{\xi})$, $1 \leq k \leq N_P$, of the expansion (\ref{eq:paper_alberto_3}) is calculated as  
\begin{multline}
\boldsymbol{\Psi}_{k}(\boldsymbol{\xi}) = \boldsymbol{\Psi}_k(\xi_1,\xi_2,\ldots,\xi_{N_U}) = \prod_{i=1}^{N_U} \psi_{\mathcal{I}_k^i}^{i}(\xi_i),
\label{eq:poly_prod}
\end{multline}
with $\sum_{i=1}^{N_U} \mathcal{I}_k^i \leq N_P$, where  $\mathcal{I}_k^i$ is an index matrix representing a graded lexicographic ordering, namely the rows of  matrix $\mathcal{I}_k^i$ 
show which order of each univariate polynomial contributes to a specific multivariate polynomial.

As a consequence of the orthonormality of the polynomials, the $N_P$ coefficients $\alpha_k(t)$ introduced in (\ref{eq:paper_alberto_3}) can be calculated as follows:
\begin{equation}
\alpha_k(t)=\int_{\boldsymbol{\xi} \in \Omega} x(t, \boldsymbol{\xi})\boldsymbol{\Psi}_k(\boldsymbol{\xi})dP(\boldsymbol{\xi}).
\label{eq:paper_alberto_4}
\end{equation}
There are different approaches to solving the integral  (\ref{eq:paper_alberto_4}), such as Galerkin projection, collocation, or numerical integration.
In particular, in this article, a Gaussian quadrature rule based on the statistical moments of $\boldsymbol{\xi}=(\xi_1, \xi_2, \ldots,\xi_{N_U})$ is used.

For a given multivariate function $\mathcal{F}(\boldsymbol{\xi})$, the full tensor product quadrature formula  $\mathcal{F}(\boldsymbol{\xi})$ can be written as
\begin{multline}
\int_{c_1}^{d_1}\cdots \int_{c_{N_U}}^{d_{N_U}} \mathcal{F}(\boldsymbol{\xi}) \approx \\
\sum_{i_1=1}^{p_1} \cdots \sum_{i_{N_U}=1}^{p_{N_U}}\mathcal{F}(\zeta^{i_1}, \ldots, \zeta^{i_{N_U}})(w_{i_1}\otimes \cdots \otimes w_{i_{N_U}}),
\label{eq:paper_alberto_5}
\end{multline}
being $\zeta^{i_j}$ and $w_{i_j}$, $1 \leq i \leq N_U$, $1\leq j \leq p$, respectively, the nodes and weights of the numerical integration, which are  obtained from the statistical moments of each random variable
$\xi_i, 1 \leq i \leq N_U$.

Let $\xi \in \Omega$ be an arbitrary component of the random vector $\boldsymbol{\xi}=(\xi_1, \xi_2, \ldots,\xi_{N_U})$, where the subindex $i$ of the random variable $\xi_i$ is omitted for the sake of clarity of exposition.
According to  \cite{oladyshkin2012data},  the  Hankel matrix of moments can be defined as
\begin{eqnarray}
M =
\begin{bmatrix}
    \mu_{0}       & \mu_{1}  & \cdots & \mu_{p} \\
    \mu_{1}       & \mu_{2}  & \cdots & \mu_{p+1} \\
    \vdots          &    \vdots &    \ddots     &    \vdots \\
    \mu_{p}       & \mu_{p+1}  & \cdots & \mu_{2p} \\
\end{bmatrix},
\label{eq:hankel_matrix}
\end{eqnarray}
where  $\mu_k$, $k \leq 0 \leq 2p$, denotes the $k$th raw statistical moment of the random variable $\xi$.

Since the Hankel matrix (\ref{eq:hankel_matrix}) is positive definite, a Cholesky decomposition 
$M = R^TR$ can be computed, such that
\begin{eqnarray*}
R =
\begin{bmatrix}
    r_{11}       & r_{12} & \cdots & r_{1,p+1} \\
                    & r_{22} & \cdots & r_{2,p+1} \\
                    &        0        &  \ddots &    \vdots \\
                   &                   &             & r_{p+1,p+1} \\
\end{bmatrix}.
\end{eqnarray*}
Let the inverse of matrix $R$ be
\begin{eqnarray*}
R^{-1} =
\begin{bmatrix}
    r^{*}_{11}      & r^{*}_{12}  & \cdots & r^{*}_{1,p+1} \\
                    & r^{*}_{22} & \cdots & r^{*}_{2,p+1} \\
                    &       0                &  \ddots &    \vdots \\
                    &                       &             & r^{*}_{p+1,p+1} \\
\end{bmatrix}.
\end{eqnarray*}
Then, each univariate orthonormal polynomial $\psi_{j}(\xi)$ used in (\ref{eq:poly_prod}) can be built from the entries of the matrix $R^{-1}$ as follows:
\begin{multline}
\psi_{j}(\xi)= r^{*}_{0j} \xi^{0} + r^{*}_{1j} \xi^{1} + r^{*}_{2j} \xi^{2} + \cdots + r^{*}_{jj} \xi^{j}, 
\label{eq:rpoly}
\end{multline}
for $0 \leq j \leq p$.

In practice, the 
coefficients $r^{*}_{ij}$ of the orthonormal polynomial (\ref{eq:rpoly}) 
can be derived directly from the elements $r_{ij}$ of the matrix $R$, avoiding the inversion of the matrix $R$.
The entries $r_{ij}$ are used to compute the coefficients $a_j$ and $b_j$ of the  three-term recurrence relation
\begin{eqnarray}
\xi\psi_{j-1} (\xi)=  b_{j-1}\psi_{j-2} (\xi) + a_{j}\psi_{j-1} (\xi) + b_{j}\psi_{j} (\xi),  
\label{eq:3term_recurrence}
\end{eqnarray}
for $1 \leq j \leq p$, where $\psi_{-1}(\xi)=0$ and $\psi_0(\xi)= 1$.
More specifically, the coefficients $a_j$ and $b_j$ are written in terms of the elements $r_{ij}$ as
 \begin{eqnarray*}
a_{j}= \frac{r_{j,j+1}}{r_{j,j}} -  \frac{r_{j-1,j}}{r_{j-1,j-1}} \quad \text{and} \quad b_{j}= \frac{r_{j+1,j+1}}{r_{j,j}},  
\end{eqnarray*}
for $1 \leq j \leq p$, with $r_{0,0}=1$ and $r_{0,1}=0$.

The coefficients $a_j$ and $b_j$ can be used to efficiently calculate the nodes and weights
used in the 
Gaussian quadrature  (\ref{eq:paper_alberto_5}).
In particular, the three-term recurrence relation (\ref{eq:3term_recurrence}) can be 
represented by means of the Jacobi matrix
\begin{eqnarray*}
J =
\begin{bmatrix}
    a_{1}      &b_{1} &  &      &   &     \\
        b_{1}              & a_{2} & b_{2} &  & 0 & \\
         &     b_{2} & a_{3} & b_{3 }&  & \\
                    &            &      \ddots       &  \ddots &    \ddots & \\
                    &      0      &          &        b_{p-2}       & a_{p-1} &  b_{p-1} \\
                                        &            &           &        &  b_{p-1}    & a_{p} \\
\end{bmatrix},
\end{eqnarray*}
which is a positive definite symmetric tri-diagonal matrix.
Then, the nodes $\zeta^i$, $1 \leq i \leq p$, associated with the univariate polynomial of order $p$ are  directly computed as the eigenvalues   of the matrix $J$, whereas the 
corresponding weights are calculated as
\begin{eqnarray*}
\omega_{i} = \nu_{1,i}^{2},
\end{eqnarray*}
where $\nu_{1,i}$, $1 \leq i \leq p$, denotes  the first component of the normalized eigenvector related to the $i$th eigenvalue of the Jacobi matrix $J$.

Therefore, the statistical moments of the random variables of the aircraft  trajectory planner model are used to determine the optimal nodes and weights of the surrogate model, as well as the corresponding orthonormal polynomials. Thus, the aPC approach offers a general framework, which permits to handle both random variables with known parametric distributions and data sets with unknown parametric distributions \cite{oladyshkin2012data}.

\section{Uncertainty Quantification}
\label{section:UQ}

The polynomial expansion (\ref{eq:paper_alberto_3}) provides a computationally efficient way to obtain the mean and variance of the output variables $x(t,\boldsymbol{\xi})$ of the aircraft trajectory planner model using the coefficients $\alpha_k(t)$. Specifically,
\begin{eqnarray*}
\mu_x &=& \alpha_1(t), \\
\sigma_x^2 &=& \sum_{k=2}^{N_P}\alpha_k^2(t).
\end{eqnarray*}
Moreover, using the quadrature rule (\ref{eq:paper_alberto_5}), the statistics of each output variable $\bold{x}(t,\boldsymbol{\xi})$ that are based on the computation of its moments, $E\left[x(t,\boldsymbol{\xi})^k\right]$, can be easily calculated in terms of scalar products, since
\begin{multline*}
E\left[x(t,\boldsymbol{\xi})^k\right] = \\
\int_{c_1}^{d_1}    \cdots \int_{c_{N_U}}^{d_{N_U}}   \left( x(t,\boldsymbol{\xi}) - E\left[  x(t,\boldsymbol{\xi})^{k-1}\right]  \right) dP(\boldsymbol{\xi}), 
\end{multline*}
for every  $k \in \mathbb{N}$. In particular, the mean and variance of each output variable $x(t,\boldsymbol{\xi})$  can also be computed as
\begin{eqnarray*}
\mu_x &=& E\left[x(t,\boldsymbol{\xi}) \right] = x(t; \zeta^{i_1}, \ldots, \zeta^{i_{N_U}}) \cdot \boldsymbol{w}, \\
\sigma_x^2 
&=& E\left[x(t,\boldsymbol{\xi})^2\right] - E\left[x(t,\boldsymbol{\xi})\right]^2 \\
&=& (x(t;\zeta^{i_1}, \ldots, \zeta^{i_{N_U}})-\mu_x)^2\cdot \boldsymbol{w},
\end{eqnarray*}

\noindent
being $\boldsymbol{w}=(w_{i_1}, \ldots, w_{i_{N_U}})$, $1 \leq   i \leq p$, the vector of Gaussian quadrature weights.

The polynomial expansion (\ref{eq:paper_alberto_3}) also allows the PDFs of the output variables of the aircraft trajectory planner model to be estimated at each instant in time. More specifically, a kernel density estimation approach can be used, which is based on the formulation of a data smoothing problem. 

Given a fixed instant in time $t^* \in [t_I,t_F]$, let 
$g_{t^*}(x)$ denote the PDF of the output variable $x(t^*,\boldsymbol{\xi}) = x(t^*; \xi_1,  \xi_2,\ldots, \xi_{N_U})$,
and let $\left\{\zeta^{i_1}, \zeta^{i_2},\ldots, \zeta^{i_{N_U}} \right\}, 1 \leq i \leq q$, be a set of samples of the vector of  random variables $\boldsymbol{\xi}=(\xi_1, \xi_2, \ldots, \xi_{N_U})$.
Then, the kernel density estimate of the PDF 
$g_{t^*}(x)$,  is calculated as
%
\begin{multline}
\label{eq:kde}
\hat{g}_{t^*}(x) = \\ \frac{1}{\eta q } \sum_{i=1}^{q} K\left( \frac{x - x(t^{*};\zeta^{i_1},\zeta^{i_2}, \ldots, \zeta^{i_{N_U}}) }{\eta} \right),
\end{multline}
where $K(\cdot)$ denotes the kernel function and $\eta$ is an appropriate smoothing parameter usually referred to as kernel bandwidth.
In particular, in this article, a Gaussian kernel is used and the corresponding kernel bandwidth is learned from the samples by means of the Silverman's rule \cite{gramacki2018nonparametric}.

Thus, given a safety threshold, the marginal and joint PDFs of the output variables can be used to estimate the probability of loss of separation between aircraft along their whole trajectories.

\section{Application}
\label{section:numerical_results}


To demonstrate the effectiveness of the probabilistic methodology for aircraft CD proposed in this paper, two numerical experiments have been conducted. In these experiments, a conflict scenario involving three aircraft has been considered, in which the trajectories of each aircraft are predicted using an aircraft trajectory planner and, given a minimum separation distance, marginal and conditional probabilities of conflict between each pair of aircraft at a given instant in time are estimated.



%
%

An aircraft trajectory planner based on optimal control has been employed to predict aircraft trajectories.
The optimal control formulation of the problem allows accurate aircraft dynamic models to be used, which is essential to improve predictability of the trajectories and obtain realistic estimates. 


In general, an Optimal Control Problem (\textsf{OCP})  \cite{benasher:2009:octwaa},\cite{betts2010practical}, is described by an objective functional that represents the  performance index of the problem, a set of Differential-Algebraic Equations (\textsf{DAE}) which specify the dynamics of the system, a set of algebraic constraints, which are usually expressed as simple bound constraints,
%
and initial and final boundary conditions. 


The objective functional $J_p$ to be minimized in the OCP formulated and solved to predict the trajectory of aircraft $p, 1 \leq p \leq N_a$, is defined as follows
\begin{equation}
J_p = \alpha_{t_p} \cdot t_{{flight}_p} + \alpha_{f_p} \cdot m_{f_p}, 
\label{eq:obj_function}
\end{equation}  
where  $N_a$ is the number of flights, $t_{{flight}_p}$ and $m_{f_p}$, are the flight time of flight $p$ and the fuel consumption of aircraft $p$, respectively, and  $\alpha_{t_p}$ and $\alpha_{f_p}$ with $\alpha_{t_p} + \alpha_{f_p} =1$ are the time and the fuel burn weighting parameters for flight $p$, respectively. These parameters are chosen based on the time and fuel-related costs incurred by airlines.

Following \cite{hull:2007:foafm}, a symmetric flight without sideslip and restricted to the horizontal plane at  cruise altitude is assumed in this paper. Therefore, a two-degree-of-freedom point variable-mass dynamic model is obtained, in which
the set of kinematics and dynamics 
\textsf{DAE} that describe the motion of the  aircraft is
\begin{eqnarray} \label{eq:kinematics_dinamics_eqs}
\nonumber  \dot{\phi}(t) &   =  &   \dfrac{V(t) \cdot \cos \chi(t) + {V_{W_N}}(t)}{R_E + h},\\  
  	\dot{\lambda}(t) &   = &   \dfrac{V(t) \cdot \sin \chi(t) + {V_{W_E}}(t)}{ \cos \phi(t) \cdot \left( R_E + h \right) }, \\  
\nonumber  \dot{\chi} (t) &   =  &  \dfrac{L(t) \cdot \sin \mu(t)}{V(t) \cdot m(t)},\\ 
\nonumber   \dot{V}(t) &   =  &   \dfrac{T(t) - D(t)  }{m(t)}, \\
\nonumber \dot{m}(t) &   =  & - T(t) \cdot  \eta(t), 
\end{eqnarray}
where the state vector has five components: the latitude $\phi$, the longitude $\lambda$, the heading angle $\chi$, the true airspeed $V$, and the mass of the aircraft $m$. In this set of equations, the control vector has three components: the thrust force $T$, the lift coefficient $C_L$, and the bank angle $\mu$. 


Other parameters and variables include the aerodynamics forces, the lift force $L$ and the drag force $D$,
the air density $\rho$, the reference wing area $S$,
the Earth radius $R_E$, the gravitational acceleration $g$, the cruise altitude $h$, and the thrust specific fuel consumption $\eta$.  Finally, ${V_{W_E}}$ and ${V_{W_N}}$ are the components of the wind velocity vector in eastward and northward directions, respectively, obtained from the 50-member \textsf{ECMWF EPS}.
The Eurocontrol's Base of Aircraft Data (\textsf{BADA}), version 3.6 \cite{bada2013}, has been used to determine the parameters of the aircraft. Specifically, the parameters of the Airbus 330-200 aircraft have been used.


The algebraic constraints considered in the problem are the flight envelope constraints which model aircraft performance limitations. These constraints are expressed in the following form
\begin{equation} \label{flight_envelope_eqs}
\eta_l(t) \leq \eta(t) \leq \eta_u(t),
\end{equation}
where $\eta(t)$ is the state or control variable and $\eta_l(t)$ and $\eta_u(t)$ are the  minimum and maximum permitted values of the variable, respectively. In particular, the upper and lower bounds for the flight envelope constraints involving the calibrated airspeed, the mass of the aircraft, the operating Mach number, the available engine thrust, the lift coefficient, and the bank angle have been applied. The explicit form of these constraints can be found in \cite{cerezo2021formation}.
Finally, the initial and final boundary conditions include the latitudes and longitudes of the initial and final positions for the three aircraft.

The described \textsf{OCP} has been transcribed into a NonLinear Programming (\textsf{NLP}) problem by discretizing the time interval $[t_I, t_F]$, $t_F = \max \{t_{{flight}_1}, \ldots, t_{{flight}_p}  \}$, and using a numerical pseudospectral method, which is described in  \cite{benson2006direct}, \cite{Garg2017AnOO}.


%
%

\begin{figure}[h!]
\center
\includegraphics[width=\columnwidth]{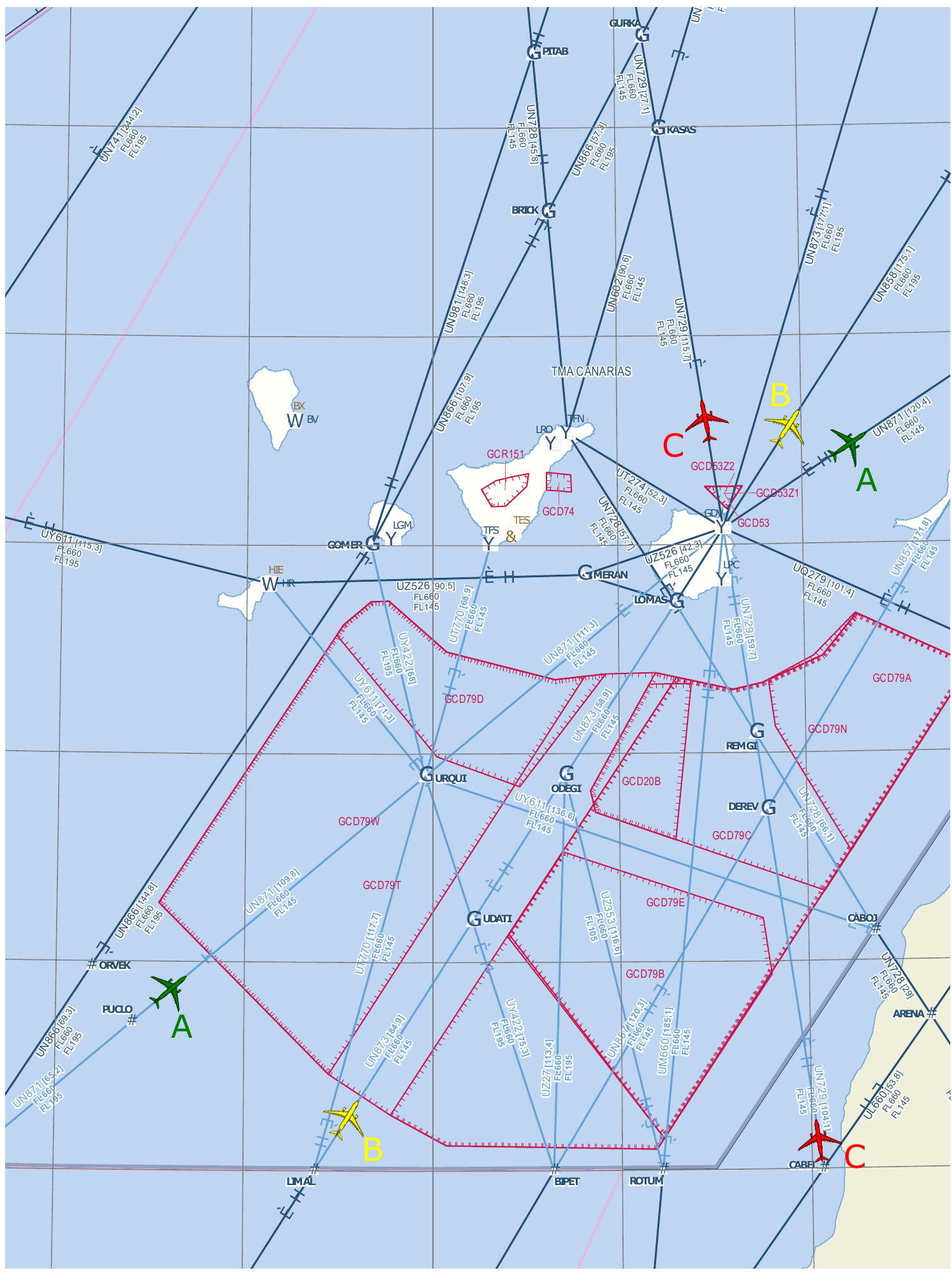}
\caption{Initial and final positions of the three aircraft considered in the numerical experiment. 
Aircraft A, B, and C are represented in green, yellow, and red colours, respectively. Source: Spanish AIP service. Not for operational use.}
\label{figure:scenario}
\end{figure}

In the considered scenario, the three aircraft,  denoted as Aircraft A, Aircraft B, and Aircraft C, are assumed to be flying at cruise level.
The latitudes and longitudes of the initial and final positions of each aircraft are given in Table~\ref{table:IFPositions}.
A representation of this scenario, along with the initial and final positions of the aircraft, is shown in Figure~\ref{figure:scenario}.
More specifically, as shown in the instrumental chart of the Spanish upper airspace represented in Figure~\ref{figure:scenario}, 
Aircraft A follows Airway {\tt UN871}, 
Aircraft B follows Airway {\tt UN873}, and, after reaching the VOR/DME {\tt GDV} of Gran Canaria, switches to Airway {\tt UN858}, and 
Aircraft C follows Airway {\tt UN729}.

\begin{table}[h!]
\center
\begin{scriptsize}
\caption{Initial and final positions of Aircraft A, Aircraft B, and Aircraft C in the numerical experiments.}
\begin{tabular}{|l c c c c|}
 \hline
  & Units & Aircraft A & Aircraft B & Aircraft C \\ [0.5ex] 
 \hline
 Initial latitude $\phi_I$ & [deg] & $25.869$ N & $25.283$ N & $25.147$ N \\ 
  \hline
 Final latitude $\phi_F$ & [deg] & $28.505$ N & $28.689$ N & $28.746$ N \\
 \hline
 Initial longitude $\lambda_I$ & [deg] & $-18.389$ E & $-17.428$ E & $-14.964$ E \\
 \hline
 Final longitude $\lambda_F$ & [deg] & $-14.677$ E & $-14.967$ E & $-15.547$ E \\
 \hline
\end{tabular}
\label{table:IFPositions}
\end{scriptsize}
\end{table}

Their trajectories are predicted using the aforementioned aircraft trajectory planner based on optimal control, in which - for all the aircraft - the time and fuel burn weighting parameters, $\alpha_{t_p}$ and $\alpha_{f_p}$, in the objective functional have been set to 1 and 0, respectively, i.e., it has been assumed that all the aircraft must reach their final positions in minimum time.
Notice that, due to the presence of wind, in general, the predicted trajectories deviate from the airways represented in Figure~\ref{figure:scenario}.


As mentioned in Section~\ref{sect:intro}, the aircraft trajectories have been predicted using wind velocity data obtained from the 50-member ECMWF EPS. 
The default eastward and northward components of the wind velocity of each member of the ECMWF EPS forecast are provided on a regular latitude-longitude grid, with a spatial resolution ranging from $0.5$ to $3$ [deg], at 9 atmospheric levels corresponding to different pressure levels ranging from $1000$ to $50$  [hPa]. Custom spatial resolutions can also be selected. The temporal coverage is 6-hour data, twice a day. The chosen pressure level is $200$ [hPa], which corresponds to the cruise altitude, and the selected spatial resolution is $0.5 \times 0.5$ [deg].

A total of 300 members of these ensembles have been employed for each numerical experiment. All of them give wind velocity forecasts for the initial time of the experiment and are valid for the duration of the experiment. Therefore, the muKL expansion of the wind velocity components has been computed using  300 realizations of the correlated random processes represented by the members of the ensembles. In particular, a truncated muKL expansion, represented by Equation~(\ref{eq:paper_cho_7}), with $M = 4$, has been considered.
The choice of the value of $M$ is a trade-off between precision and computational cost.
Since, for a given value of $M$, $2^M$ aircraft trajectories must be computed, and $2^8 = 256$, a number of the same order of magnitude as the number of members of the ensembles that have been considered, computational advantages are obtained from the muKL expansion only for $M < 8$. However, for precision reasons, taking $M \leq 2$ is not adequate.

The only difference between the two numerical experiments, which will be referred to as Experiment 1 and Experiment 2, is the time at which the three aircraft are located at their initial positions, which is the initial time $t_I$ of the experiment. Therefore, in the two experiments, different wind velocity forecasts have been used.

\begin{itemize}

\item
Experiment 1: Collision detection for March 18, 2022, at 06:00 UTC.

\item
Experiment 2: Collision detection for March 19, 2022, at 06:00 UTC.

\end{itemize}

As already mentioned in Section~\ref{section:muKL}, in order to include the wind velocity information in the aircraft trajectory planner model, an analytic function that interpolates the eigenfunctions of the wind velocity components derived from the muKL expansion must be determined. More specifically, in this paper, a Radial Basis Functions (RBF)-based interpolation method is employed. The RBF is a widespread technique that allows multidimensional data to be interpolated, which can be gridded or scattered, i.e. the RBF-based method can interpolate both structured and unstructured data \cite{buhmann2003radial}.
In the numerical experiments, a series of equally weighted forecasts with different time horizons have been used.


Once the aircraft trajectories are predicted for each instant of the time discretization
using the optimal control-based aircraft trajectory planner, the corresponding optimal state variables are used to calculate the main statistics of the distance between each pair of aircraft at each instant of the time discretization. More specifically, the orthodromic distance between aircraft is computed from the obtained optimal latitudes and longitudes using the haversine formula. Then, the 2-sigma confidence envelopes are calculated. If the 2-sigma confidence envelope  of the distance between two aircraft crosses the threshold at a given instant in time, a conflict is detected and no further analysis is carried out. In particular, in this case, the probability of conflict is not computed because the value of this probability would only confirm the existence of a conflict.  
In all other cases, the PDF of the distance between each pair of aircraft 
at that instant in time is determined, from which the probability of conflict is calculated. 
In some cases, this analysis confirms the absence of conflict, whereas in other cases,  
a conflict is detected.

\subsection{Experiment 1: Collision detection for March 18, 2022 at 06:00 UTC}

\begin{figure}[h!]
\psfrag{K}{\small \sf $mg$}
\begin{subfigmatrix}{2}
\subfigure[Eastward component.]{\includegraphics[width=0.48\columnwidth]{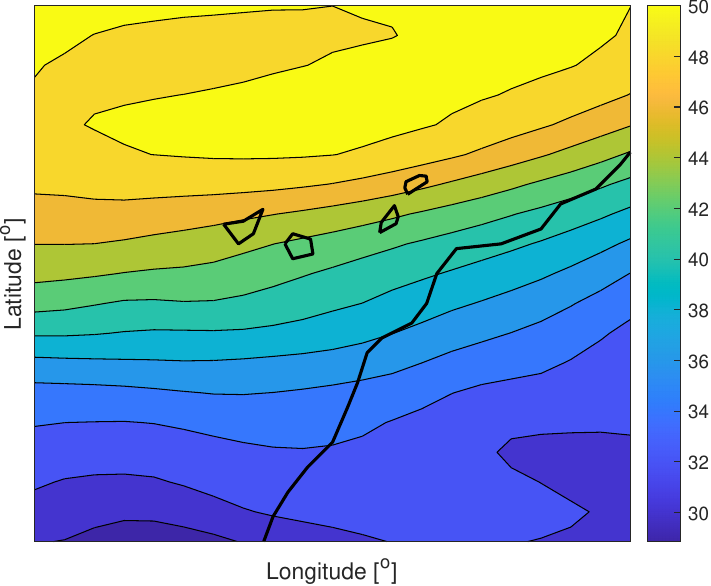}}\hspace{0 cm}
\subfigure[Northward component.]{\includegraphics[width=0.48\columnwidth]{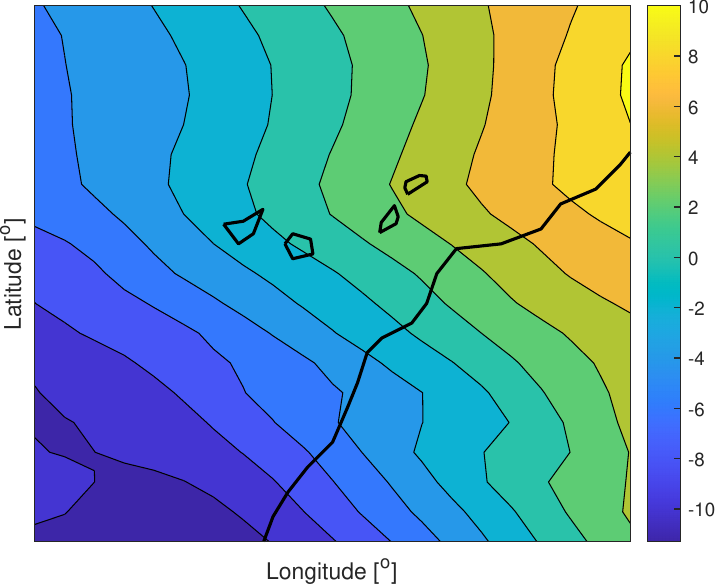}}
\end{subfigmatrix}
\caption{Experiment 1: Mean value,  in [m/s], of the magnitude of the components of the wind velocity calculated using ECMWF EPS forecasts for March 18, 2022, at 6:00 for the pressure altitude of 200 [hPa].}
\label{figure:mean_realizations_experiment_1} 
\end{figure}

\begin{figure}[h!]
\psfrag{K}{\small \sf $mg$}
 \begin{subfigmatrix}{2}
   \subfigure[Eastward component.]{\includegraphics[width=0.48\columnwidth]{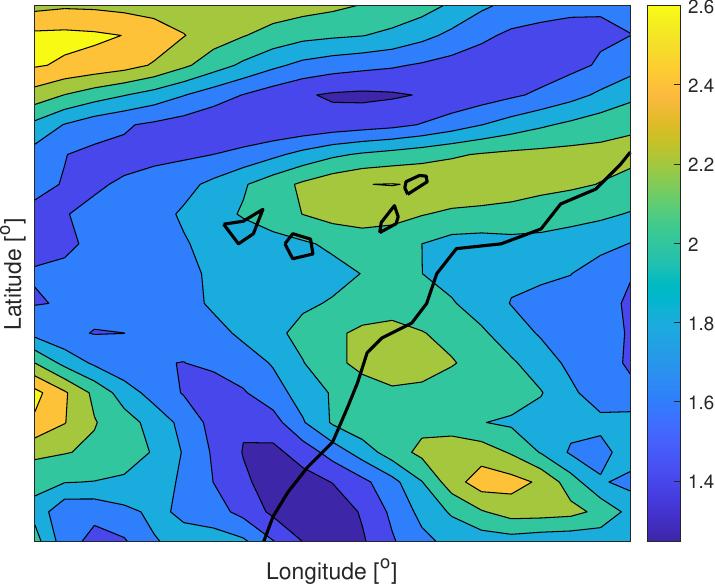}}
  \hspace{0 cm}
   \subfigure[Northward component.]{\includegraphics[width=0.48\columnwidth]{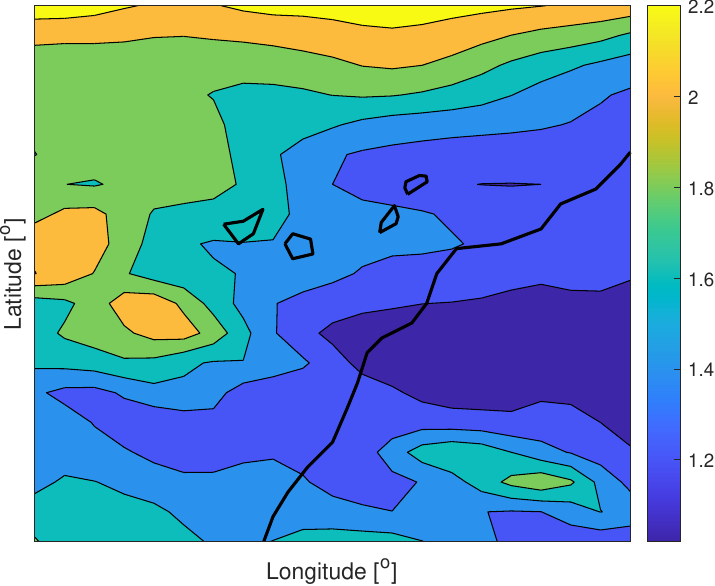}}
 \end{subfigmatrix}
 \caption{Experiment 1: Standard deviation, in [m/s], of the magnitude of the components of the wind velocity calculated using ECMWF EPS forecasts for March 18, 2022, at 6:00 for the pressure altitude of 200 [hPa].}
\label{figure:std_realizations_experiment_1} 
\end{figure}

In this experiment, the initial time is March 18, 2022, at 06:00 UTC,
and therefore, the wind velocity predictions published 
on March 18, 2022, at 00:00 with a time horizon of 06 [h],
on March 17, 2022, at 12:00 with a time horizon of 18 [h],
on March 17, 2022, at 00:00 with a time horizon of 30 [h],
on March 16, 2022, at 12:00 with a time horizon of 42 [h],
on March 16, 2022, at 10:00 with a time horizon of 54 [h], and 
on March 15, 2022, at 12:00 with a time horizon of 66 [h], have been used.
The mean and standard deviation, measured in [m/s], of the eastward and northward components of the wind velocity obtained from these 300 ensemble members are shown in  figures~\ref{figure:mean_realizations_experiment_1} and \ref{figure:std_realizations_experiment_1}, respectively.
The percentages of variability explained by each random variable 
in the truncated muKL expansion are shown in Table~\ref{table:explanation_experiment_1}.

\begin{figure}[h!]
\psfrag{K}{\small \sf $mg$}
\begin{subfigmatrix}{2}
\subfigure[Eastward component.]{\includegraphics[width=0.48\columnwidth]{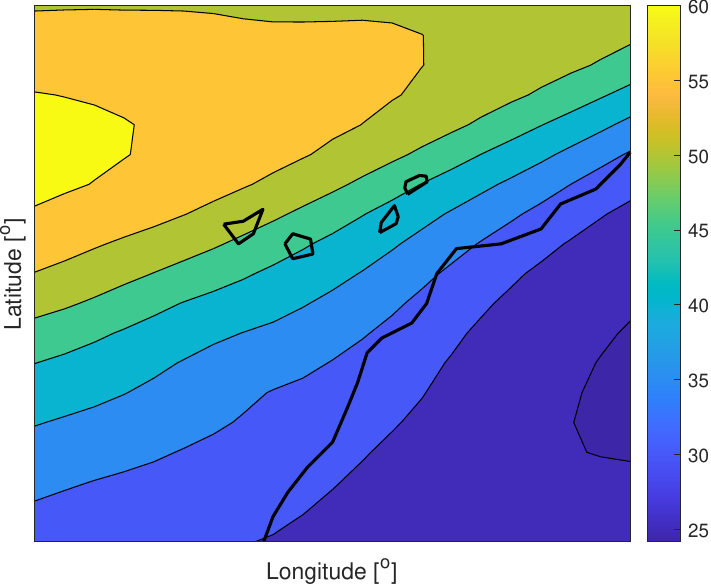}}\hspace{0 cm}
\subfigure[Northward component.]{\includegraphics[width=0.48\columnwidth]{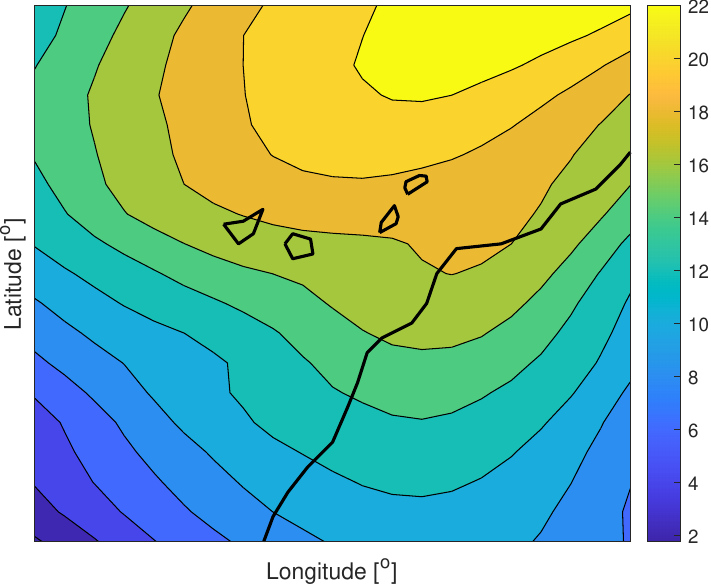}}
\end{subfigmatrix}
\caption{Experiment 2: Mean value,  in [m/s], of the magnitude of the components of the wind velocity calculated using ECMWF EPS forecasts for March 19, 2022, at 6:00 for the pressure altitude of 200 [hPa]. }
\label{figure:mean_realizations_experiment_2} 
\end{figure}

\begin{figure}[h!]
\psfrag{K}{\small \sf $mg$}
 \begin{subfigmatrix}{2}
   \subfigure[Eastward component.]{\includegraphics[width=0.48\columnwidth]{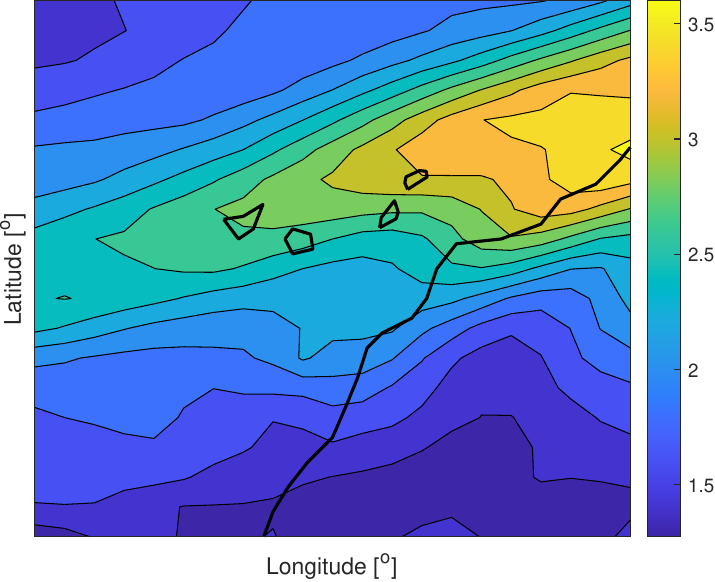}}
  \hspace{0 cm}
   \subfigure[Northward component.]{\includegraphics[width=0.48\columnwidth]{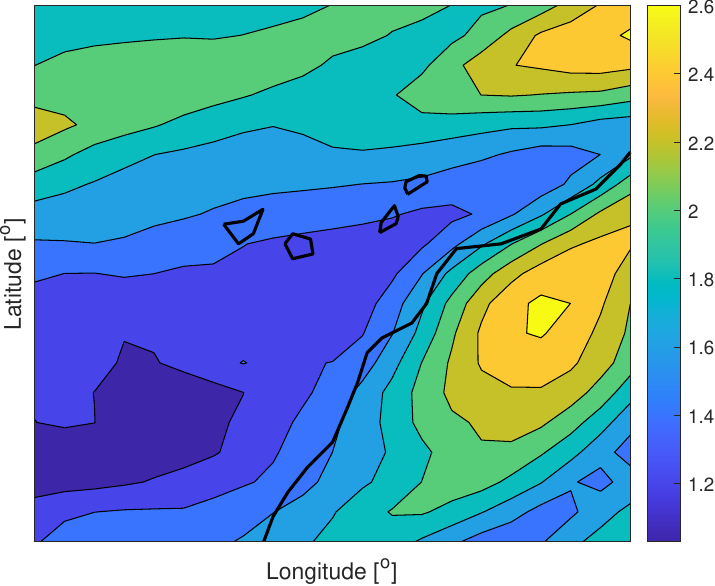}}
 \end{subfigmatrix}
 \caption{Experiment 2: Standard deviation, in [m/s], of the magnitude of the components of the wind velocity calculated using ECMWF EPS forecasts for March 19, 2022, at 6:00 for the pressure altitude of 200 [hPa].}
\label{figure:std_realizations_experiment_2} 
\end{figure}

\begin{table}
\center
\caption{Experiment 1: Percentages of explanation of each random variable of the muKL expansion.}
\begin{tabular}{| c c|}
 \hline
Random variable & \% of explanation \\  
 \hline
 $\xi_1(\omega)$  & 47.8650   \\ 
  \hline
 $\xi_2(\omega)$  & 15.6234  \\
 \hline
 $\xi_3(\omega)$  & 5. 8512   \\
 \hline
 $\xi_4(\omega)$  & 4.6298   \\
 \hline
 Total & 73.9694 \\
 \hline
\end{tabular}
\label{table:explanation_experiment_1}
\end{table}

%
%
%
%
%
%
%
%
%
%
%
%
%

\begin{figure}[h!]
\center
\includegraphics[width=\columnwidth]{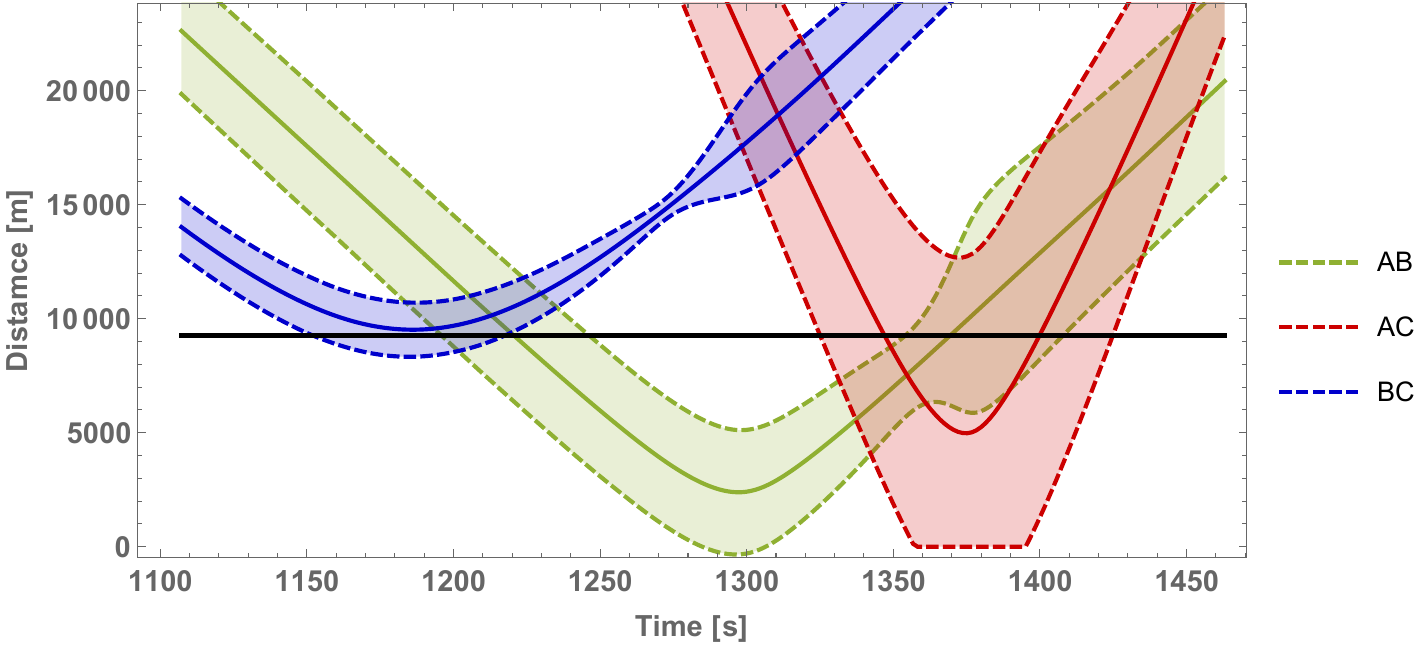}
\caption{Experiment 1: Mean distance between each pair of aircraft
along with the corresponding
2-sigma confidence envelopes. The horizontal solid black line represents a distance of 5 [NM],  the minimum separation required. }
\label{figure:confidence_envelopes_distances_experiment_1}
\end{figure}

\begin{figure}[h!]
\center
\includegraphics[width=\columnwidth]{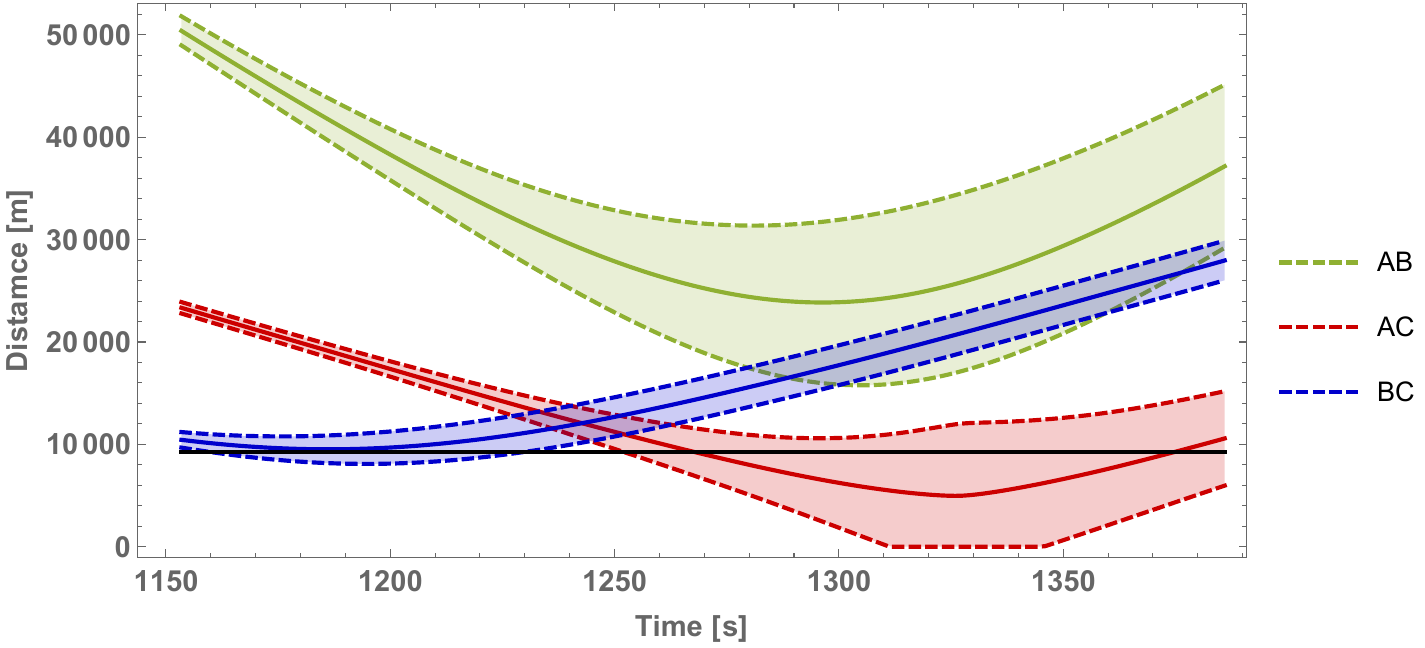}
\caption{Experiment 2: Mean distance between each pair of aircraft
along with the corresponding
2-sigma confidence envelopes. The horizontal solid black line represents a distance of 5 [NM],  the minimum separation required. }
\label{figure:confidence_envelopes_distances_experiment_2}
\end{figure}

Figure~\ref{figure:confidence_envelopes_distances_experiment_1} shows the  mean distances between each pair of aircraft, along with their associated 2-sigma confidence envelopes calculated using the proposed  CD  probabilistic methodology. The horizontal solid black line represents the minimum separation distance between aircraft required by the current regulation, namely 5 [NM] \cite{icao4444doc}.
Figure~\ref{figure:confidence_envelopes_distances_experiment_1} shows that, based on the 2-sigma confidence envelope criterion, conflicts between 
Aircraft A and Aircraft B, 
Aircraft A and Aircraft C, and  
Aircraft B and Aircraft C are detected. 
As said before, in this case, the probabilities of conflict are not determined, because their values would only confirm the existence of these conflicts.

\subsection{Experiment 2: March 19, 2022 at 06:00 UTC}

In this experiment, the initial time is March 19, 2022, at 06:00 UTC, and therefore, the wind velocity predictions published 
on March 19, 2022, at 00:00 with a time horizon of 06 [h],
on March 18, 2022, at 12:00 with a time horizon of 18 [h],
on March 18, 2022, at 00:00 with a time horizon of 30 [h],
on March 17, 2022, at 12:00 with a time horizon of 42 [h],
on March 17, 2022, at 10:00 with a time horizon of 54 [h], and 
on March 16, 2022, at 12:00 with a time horizon of 66 [h], have been used.
The mean and standard deviation, measured in [m/s], of the eastward and northward components of the wind velocity obtained from these 300 ensemble members are shown in  figures~\ref{figure:mean_realizations_experiment_2} and \ref{figure:std_realizations_experiment_2}, respectively.
The percentages of variability explained by each random variable 
in the truncated muKL expansion are shown in Table~\ref{table:explanation_experiment_2}. 

\begin{table}
\center
\caption{Experiment 2: Percentages of explanation of each random variable of the muKL expansion.}
\begin{tabular}{| c c|}
 \hline
Random variable & \% of explanation \\  
 \hline
 $\xi_1(\omega)$  & 29.394 \\ 
  \hline
 $\xi_2(\omega)$  & 17.930  \\
 \hline
 $\xi_3(\omega)$  & 13.079 \\
 \hline
 $\xi_4(\omega)$  & 8.610\\
 \hline
 Total & 69.013 \\
 \hline
\end{tabular}
\label{table:explanation_experiment_2}
\end{table}

Figure~\ref{figure:confidence_envelopes_distances_experiment_2} shows the  mean distances between each pair of aircraft, along with their associated 2-sigma confidence envelopes, calculated using the proposed  CD  probabilistic methodology. The horizontal solid black line represents the minimum separation distance between aircraft required by the current regulation, namely 5 [NM] \cite{icao4444doc}.
Figure~\ref{figure:confidence_envelopes_distances_experiment_2} shows that, based on the 2-sigma confidence envelope criterion, conflicts between Aircraft A and Aircraft C and between Aircraft B and Aircraft C are detected, whereas no conflict between Aircraft A and Aircraft B is indicated. 
Moreover, the same conclusion is obtained if a 3-sigma confidence envelope criterion is used.

In this case, since the 2-sigma confidence criterion does not detect any conflict between Aircraft A and Aircraft B, a further analysis is carried out in the vicinity of the instant in time at which the aircraft are at the minimum distance. First, the PDF of the distance 
between Aircraft A and Aircraft B are calculated, and then the probability of conflict is estimated.
\begin{figure}[h!]
\center
\includegraphics[width=\columnwidth]{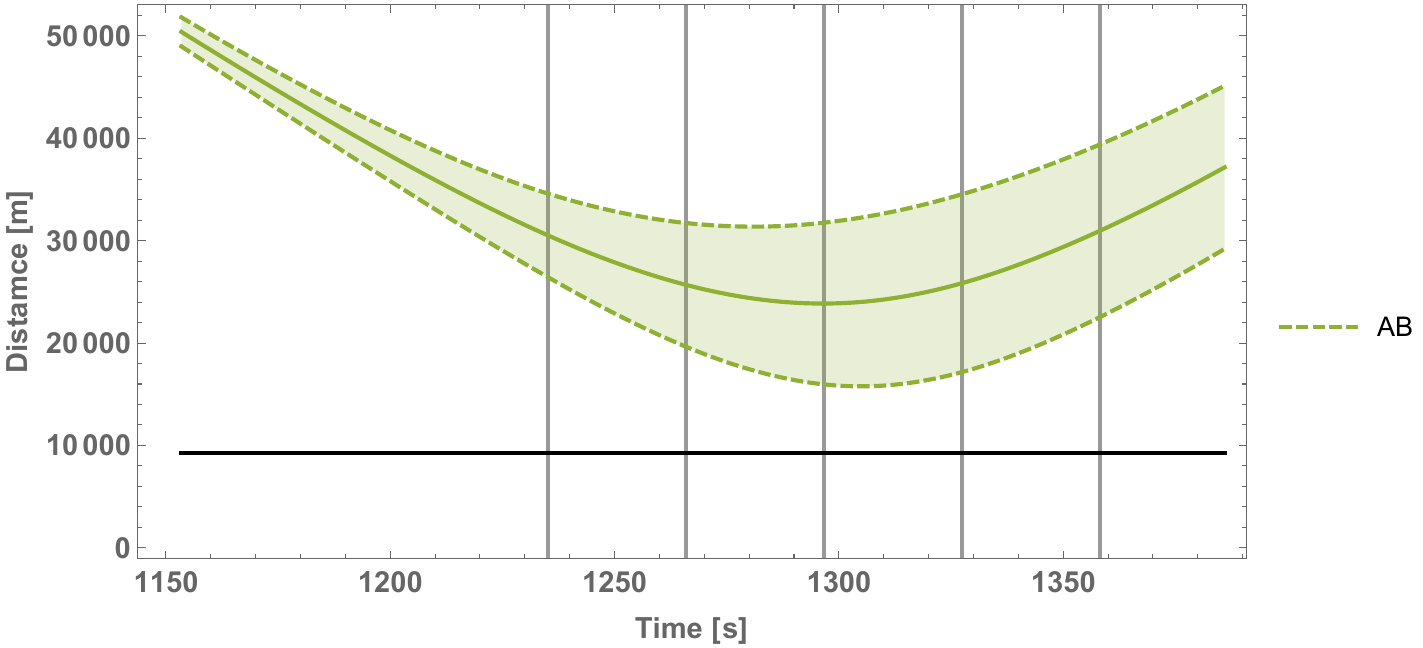}
\caption{Experiment 2: Mean distance between Aircraft A and Aircraft B, 
along with the  
corresponding 2-sigma confidence envelope. The vertical segments in grey represent five instants in time, namely $1235.08$ [s], $1265.88$ [s], $1296.68 $ [s], $1327.48$ [s], and $1358.28$ [s], at which the PDF of the distance between these aircraft is calculated and shown in Figure~\ref{figure:pdfs_distances_aircraft_A_aircraft_B_experiment_2}. The horizontal solid black line represents a distance of 5 [NM], the minimum separation required.}
\label{figure:confidence_envelope_distance_aircraft_A_aircraft_B_experiment_2}
\end{figure}
\begin{figure}[h!]
\center
\includegraphics[width=\columnwidth]{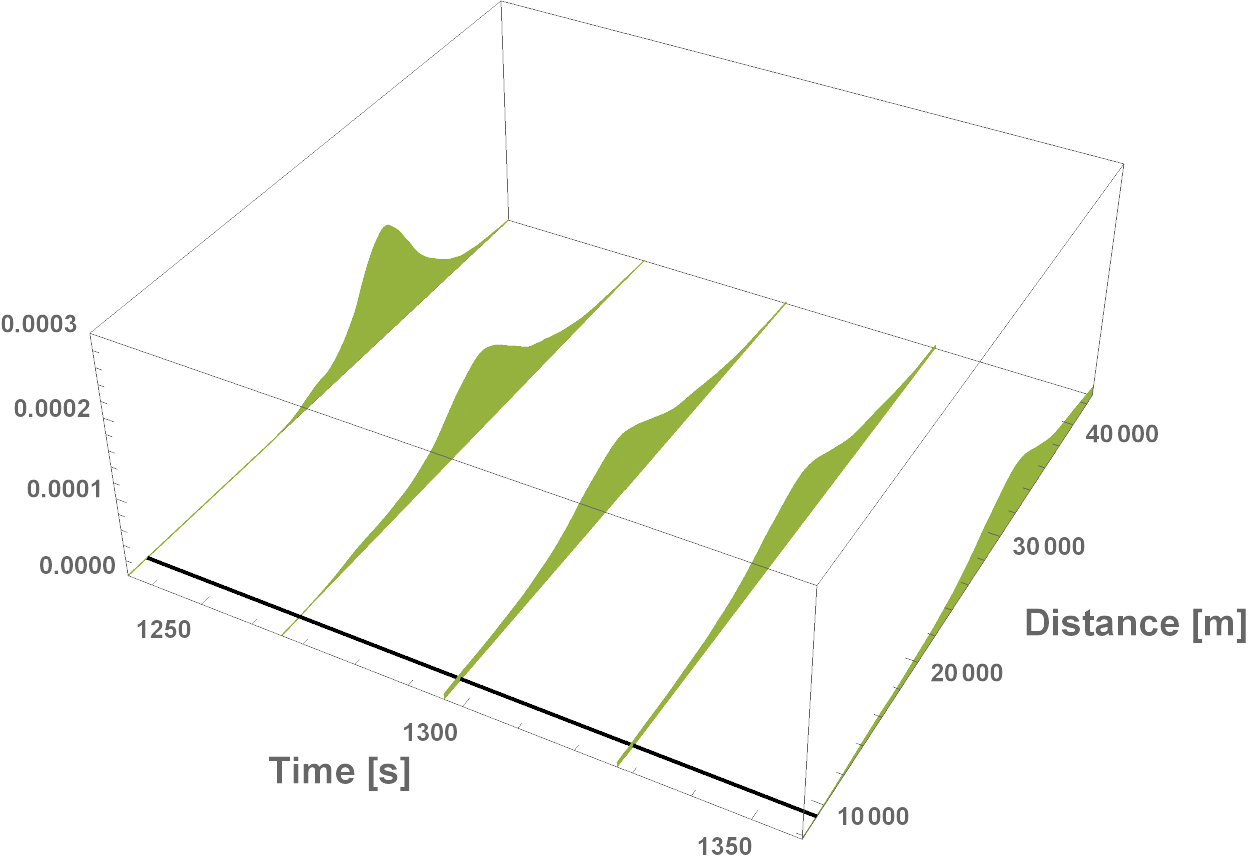}
\caption{Experiment 2: PDFs of the distance between Aircraft A and Aircraft B at instants $1235.08$ [s], $1265.88$ [s], $1296.68 $ [s], $1327.48$ [s], and $1358.28$ [s].
The solid black line represents a distance of 5 [NM], the minimum separation required.}
\label{figure:pdfs_distances_aircraft_A_aircraft_B_experiment_2}
\end{figure}
The vertical grey segments depicted in Figure~\ref{figure:confidence_envelope_distance_aircraft_A_aircraft_B_experiment_2} represent five instants in time, namely $1235.08$ [s], $1265.88$ [s], $1296.68$ [s], $1327.48$ [s], and $1358.28$ [s], at which the PDFs of the distances between Aircraft A and Aircraft B have been estimated.
The resulting PDFs are shown in Figure~\ref{figure:pdfs_distances_aircraft_A_aircraft_B_experiment_2}.
When the aircraft are at the minimum distance, which occurs at time $1296.68$ [s],  the probability of conflict is $0.0260732$, which is a significant value from the point of view of air traffic safety,
where probabilities of conflict greater than $10^{-2}$ are considered high risk of conflict.
Thus, it can be concluded that, unlike the CD criterion based on the 2-sigma confidence envelope, the CD criterion based on probability is able to detect a conflict between Aircraft A and Aircraft B.

\begin{figure}[h!]
\center
\psfrag{K}{\small \sf $mg$}
 \begin{subfigmatrix}{1}
   \subfigure[Joint PDF.]{\includegraphics[width=\columnwidth]{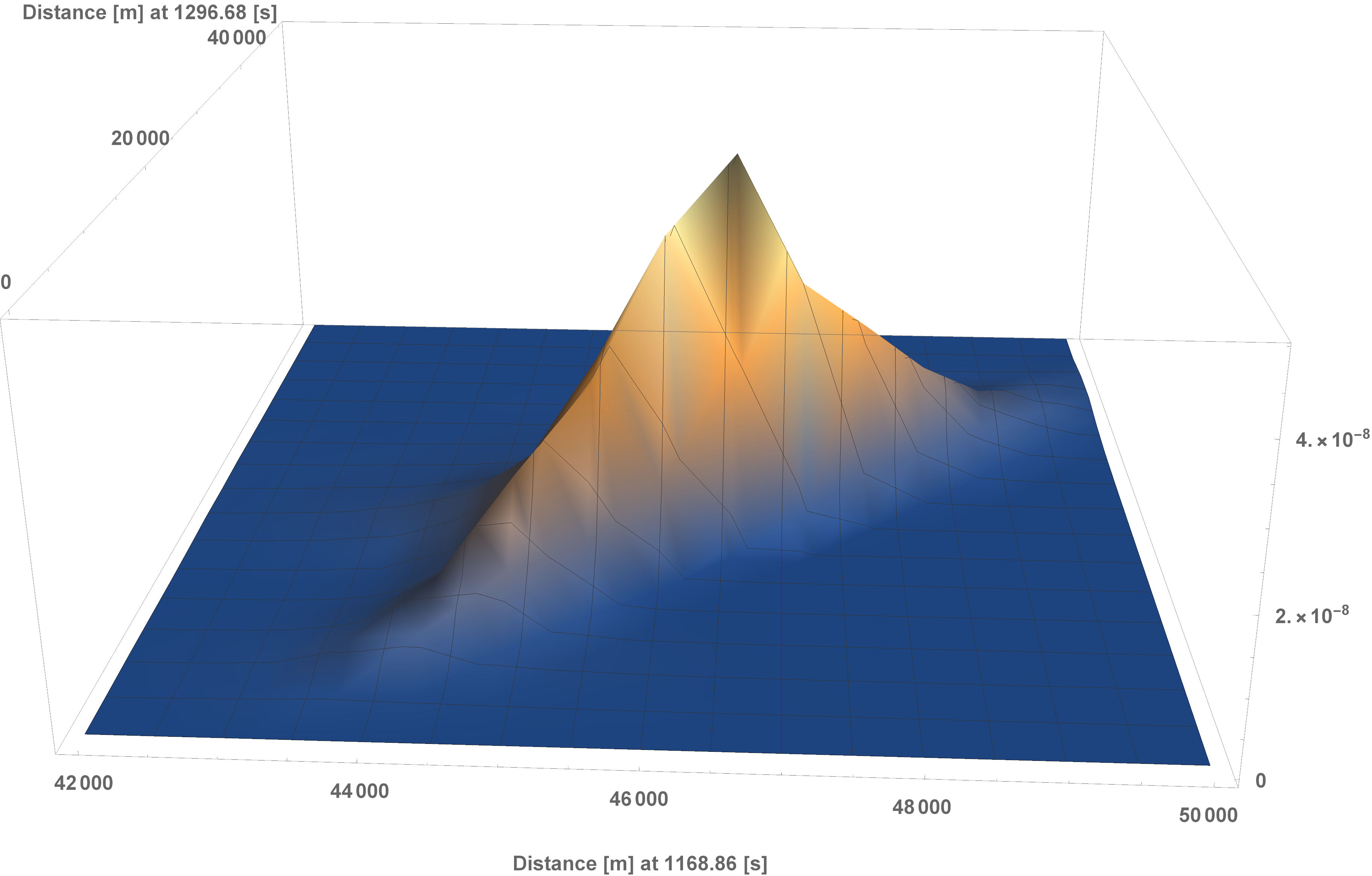}}
  \hspace{0 cm}
   \subfigure[Joint CDF.]{\includegraphics[width=\columnwidth]{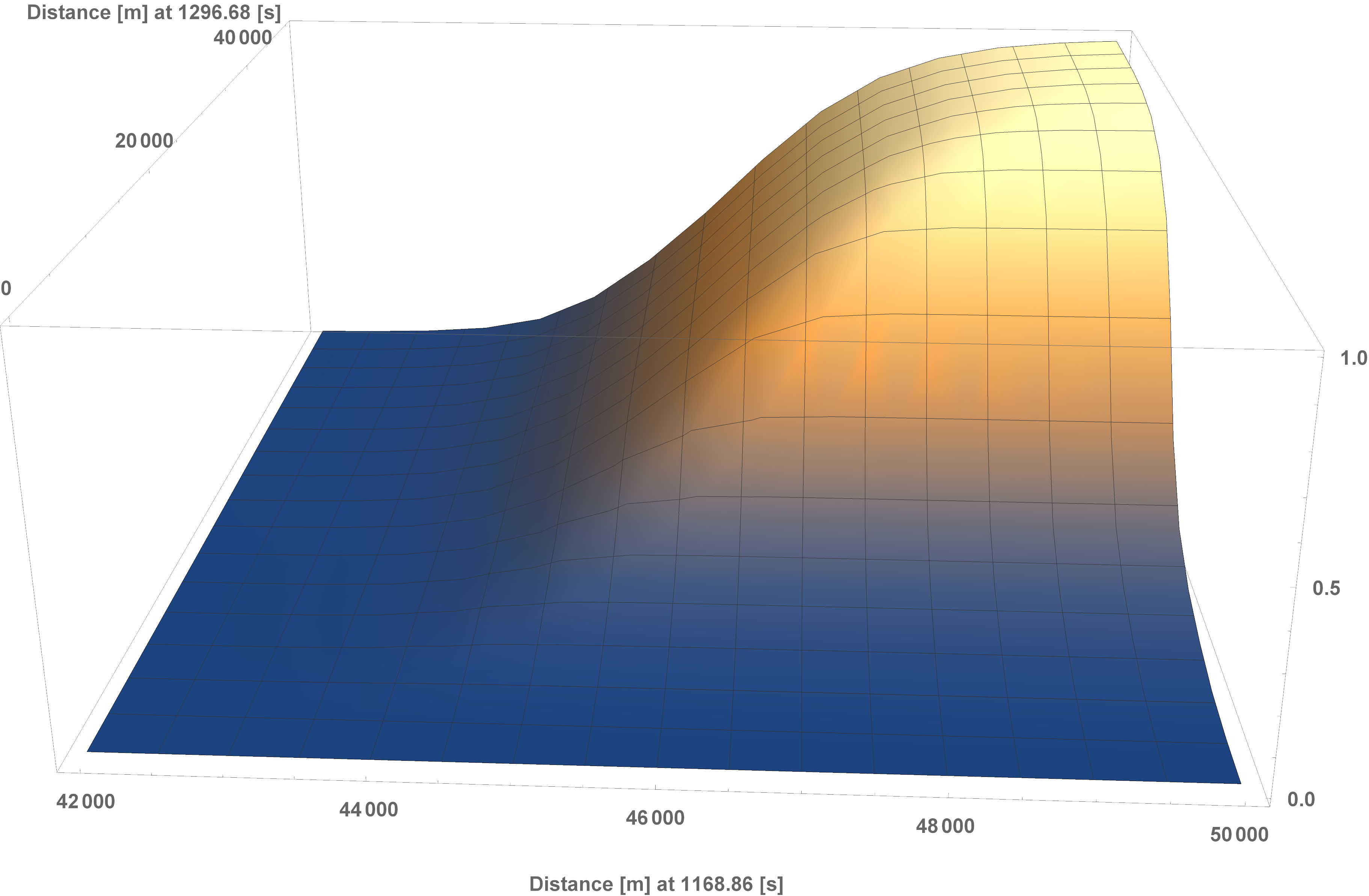}}
 \end{subfigmatrix}
 \caption{Experiment 2: Joint PDF and joint CDF of the distances between Aircraft A and Aircraft B at instants $1168.86$ and $1296.68$ [s].}
 \label{figure:bivariate_pdf_experiment_2} 
\end{figure}

The surrogate model (\ref{eq:paper_alberto_3}) allows not only the marginal PDF of the distance between two aircraft at any instant to be calculated but also the joint PDF of the distances between two aircraft at two different instants in time to be computed.  For example, the joint PDF of the distances between Aircraft A and Aircraft B at times  $1168.86$ [s] and $1296.68$ [s]  is represented in Figure~\ref{figure:bivariate_pdf_experiment_2}.a, whereas the corresponding joint CDF is shown in Figure~\ref{figure:bivariate_pdf_experiment_2}.b. From this joint probability distribution, it is possible to calculate both the joint and conditional probabilities of conflict.
For instance, knowing that the distance between Aircraft A and Aircraft B at time $1168.86$ [s] is  lower than $25$ [NM], the probability of conflict at time 1296.68 [s] conditioned by this event becomes 
$0.0496445$, which is considerably higher than  the marginal probability of $0.0260732$.

%
%

To study the effect of the truncation order of the muKL expansion, represented by Equation~(\ref{eq:paper_cho_7}),
the probability of conflict between Aircraft A and Aircraft B when the aircraft are at the minimum distance, which occurs at time $1296.68$ [s], 
has also been calculated for $M = 3$ and $M = 5$. The results are shown in the first column of Table~\ref{table:probabilities_different_M}. 
To study the effect of the truncation order of the muKL expansion, represented by Equation~(\ref{eq:paper_cho_7}),
the probability of conflict at time $1296.68$ [s], conditioned by the event that the distance between Aircraft A and Aircraft B at time $1168.86$ [s] is lower than $25$ [NM], has also been calculated for $M = 3$ and $M = 5$. The results are shown in the second column of Table~\ref{table:probabilities_different_M}.  
The mean computation times for computing the marginal and the conditional probabilities of conflict for $M = 3, 4, 5$, 
obtained by averaging over $10$ executions the computation times observed on a standard laptop computer, are reported in the third column of Table~IV. 
Notice that that the computation times for calculating both the marginal and conditional probabilities of conflict are nearly identical. Therefore, only a single mean computation time for each value of $M$ has been reported in Table~IV. 
For $M = 8$, the probability of conflict is $0.0279676$ and the conditional probability of conflict is $0.0548837$. 
As expected, the computation time for this value of $M$ increases to over $4$ hours.


\begin{table*}
\center
\caption{Experiment 2: Probabilities of conflict between Aircraft A and Aircraft B at instant $1296.68$ [s] as functions of the values of $M$ in the muKL expansion.}
\begin{tabular}{| c  c  c  c | }
\hline
 $M$ & Probability of conflict& Conditional probability of conflict & Mean computation time [s]\\ 
\hline
3 & 0.0153225 &  0.0318974  & 120.95 \\
\hline
4 & 0.0260732 &  0.0496445  & 327.29 \\ 
\hline
5 & 0.0273232 &  0.0545471  &  831.55 \\
\hline
\end{tabular}
 \label{table:probabilities_different_M}
\end{table*}

%
%

\begin{figure}[h!]
\center
\includegraphics[width=\columnwidth]{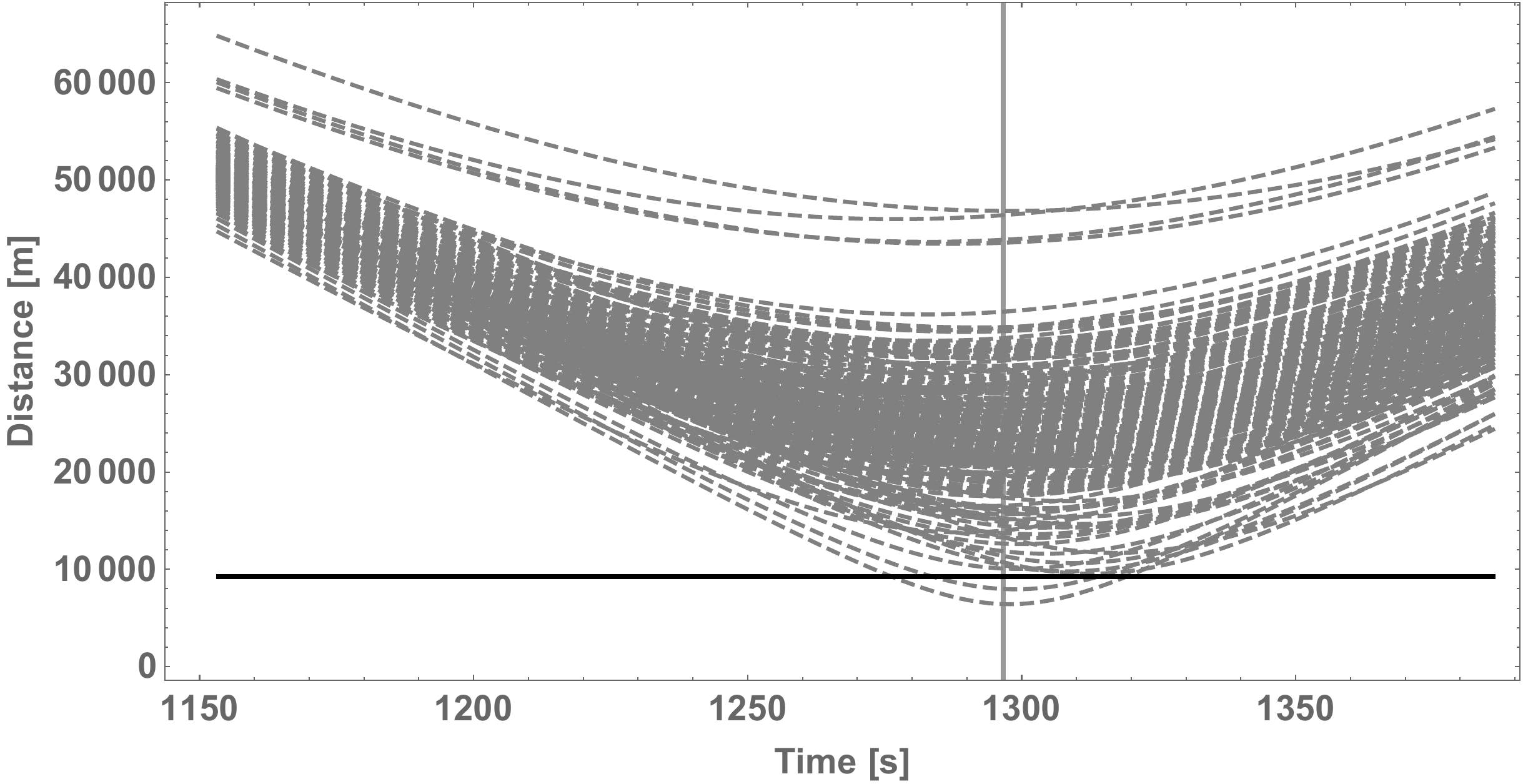}
\caption{Experiment 2: Distances between Aircraft A and Aircraft B as functions of time calculated for the $300$ pairs of trajectories generated using the ensemble approach to estimate the probability of conflict between Aircraft A and Aircraft B at instant $1296.68$ [s]. The vertical grey segment represents 
instant $1296.68$ [s]. 
The horizontal solid black line represents in [m] the distance of 5 [NM], the minimum separation required.
}
\label{figure:distances_aircraft_A_aircraft_B_ensemble_approach}
\end{figure}

The obtained results have been compared with those obtained using the ensemble method for aircraft CD in the presence of uncertainty in the wind velocity forecast \cite{eulaliaTesis}, \cite{gonzalez2018robust} described in the introduction. 
To estimate the probability of conflict between Aircraft A and Aircraft B at time $1296.68$ [s], 300 pairs of trajectories of Aircraft A and Aircraft B have been calculated deterministically, each of them using a different member of the EPS forecast. Then, the distances between Aircraft A and Aircraft B have been computed as functions of time, which are represented in Figure~\ref{figure:distances_aircraft_A_aircraft_B_ensemble_approach}.
In this figure, the vertical grey segment represents 
instant $1296.68$ [s], and
the horizontal solid black line represents a distance of 5 [NM], the minimum separation required.
The figure shows that the minimum distance between these aircraft is below the threshold of 5 [NM] in only 2 pairs of trajectories out of $300$. This corresponds to a probability of conflict between Aircraft A and Aircraft B of $0.00667$.  
The probability of conflict at time $1296.68$ [s], conditioned by the event that the distance between Aircraft A and Aircraft B at time $1168.86$ [s] is lower than $25$ [NM], has also been calculated using the ensemble method. The result is $0.0123457$, which is considerably higher than the marginal probability of $0.00667$. As expected, the ensemble method also detects a conflict between Aircraft A and Aircraft B at time $1296.68$ [s], although the marginal and conditional probabilities of conflict calculated with this method differ from those calculated with the method proposed in this paper. In particular, lower values of the probability of conflict are obtained with the ensemble method. These differences are due to the fact that, in the method proposed in this paper, the probabilities of conflict are calculated using probability distributions of distances between aircraft whereas, in the ensemble method, they are calculated based on the occurrence of events involving theses distances. It is easy to see that, in the ensemble method, information about proximity between aircraft, which is contained in the probability distributions of distances between aircraft, is not considered. This means that the proposed method is able to estimate the probability of conflict more accurately than the ensemble method, which tends to underestimate the probability of conflict.

\section{Conclusions}
\label{section:conclusions}


In this paper, a probabilistic approach to aircraft conflict detection in the cruise phase of flight in the presence of wind velocity prediction uncertainty quantified by an ensemble weather forecasts has been proposed.
Specifically, a conflict among aircraft has been assumed to occur when the distance between two or more aircraft is smaller than the required minimum separation.


For this purpose, the multiple uncorrelated Karhunen-Lo\`eve expansion, the arbitrary polynomial chaos expansion, and a Gaussian kernel density estimator have been combined with a deterministic trajectory planner to 
quantify the 
uncertainty associated with the northward and eastward components of the wind velocity
and obtain the probability density function of the distance between aircraft at different instants in time, which allows the probability of conflict between pairs of aircraft to be calculated. The joint probability density function of the distance between two aircraft at two different instants in time has also been estimated, which permits joint and conditional probabilities of conflict between pairs of aircraft to be computed.


The ensemble prediction system provided by the European Centre for Medium-Range Weather Forecast has been used in this paper as the probabilistic wind velocity forecast, whereas a deterministic aircraft trajectory planner based on pseudospectral optimal control has been used to generate the aircraft trajectories.
However, the proposed data-driven probabilistic methodology for aircraft conflict detection is a general approach in which 
any ensemble prediction system and any deterministic aircraft trajectory planner can be used.


In the numerical experiments, the proposed methodology has been applied to detect conflicts among aircraft flying at cruise altitude. 
More specifically, conflicts between aircraft have been detected not only by calculating the marginal probability of occurrence at a certain instant in time but also by computing the conditional probability of occurrence, in which the conditioning event is assumed to be information about the distance between the two aircraft at a previous instant. 
As expected, information about the distance between two aircraft at a certain previous instant strongly influences the probability of conflict between the same aircraft at a future instant.
Moreover, it has been proven that this approach to conflict detection based on the probability of occurrence is much more accurate than the approach based on confidence envelopes of the separation distance between aircraft, which is very significant in terms of safety for the air traffic management system.


The proposed methodology provides a new framework for expanding the capabilities of current aircraft conflict detection systems in the presence of wind velocity forecasting uncertainty, 
insofar as it yields marginal and conditional probabilities of conflict and not just confidence intervals.


%
%
%
%
%

\bibliographystyle{IEEEtran}
\bibliography{conflict_detection_bibliography}

\begin{IEEEbiography}[{\includegraphics[width=1in,height=1.25in,clip,keepaspectratio]{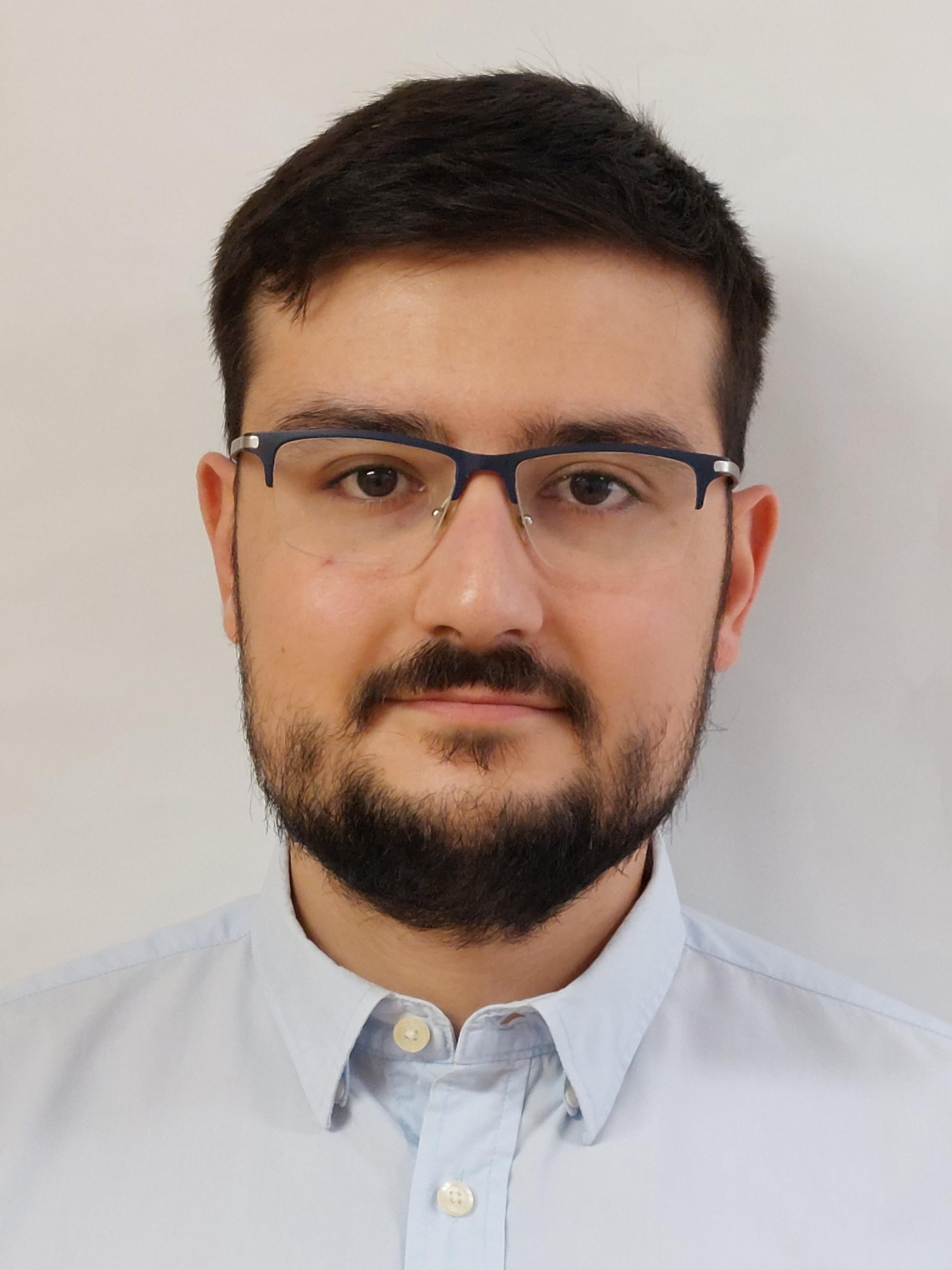}}]{Jaime de la Mota}
is a Ph.D. student at the Universidad Rey Juan Carlos, Madrid, Spain. 
He received the B.Sc. degree in physics from the Universitat de Val\`encia, Valencia, Spain, and the M.Sc. degree in meteorology and geophysics from the Universidad Complutense de Madrid, Madrid, Spain, in 2016 and 2017, respectively. 
He was a Teaching Assistant of Statistics at the Universidad Rey Juan Carlos.
His research is focused on uncertainty quantification applied to air traffic management.
\end{IEEEbiography}

\begin{IEEEbiography}[{\includegraphics[width=1in,height=1.25in,clip,keepaspectratio]{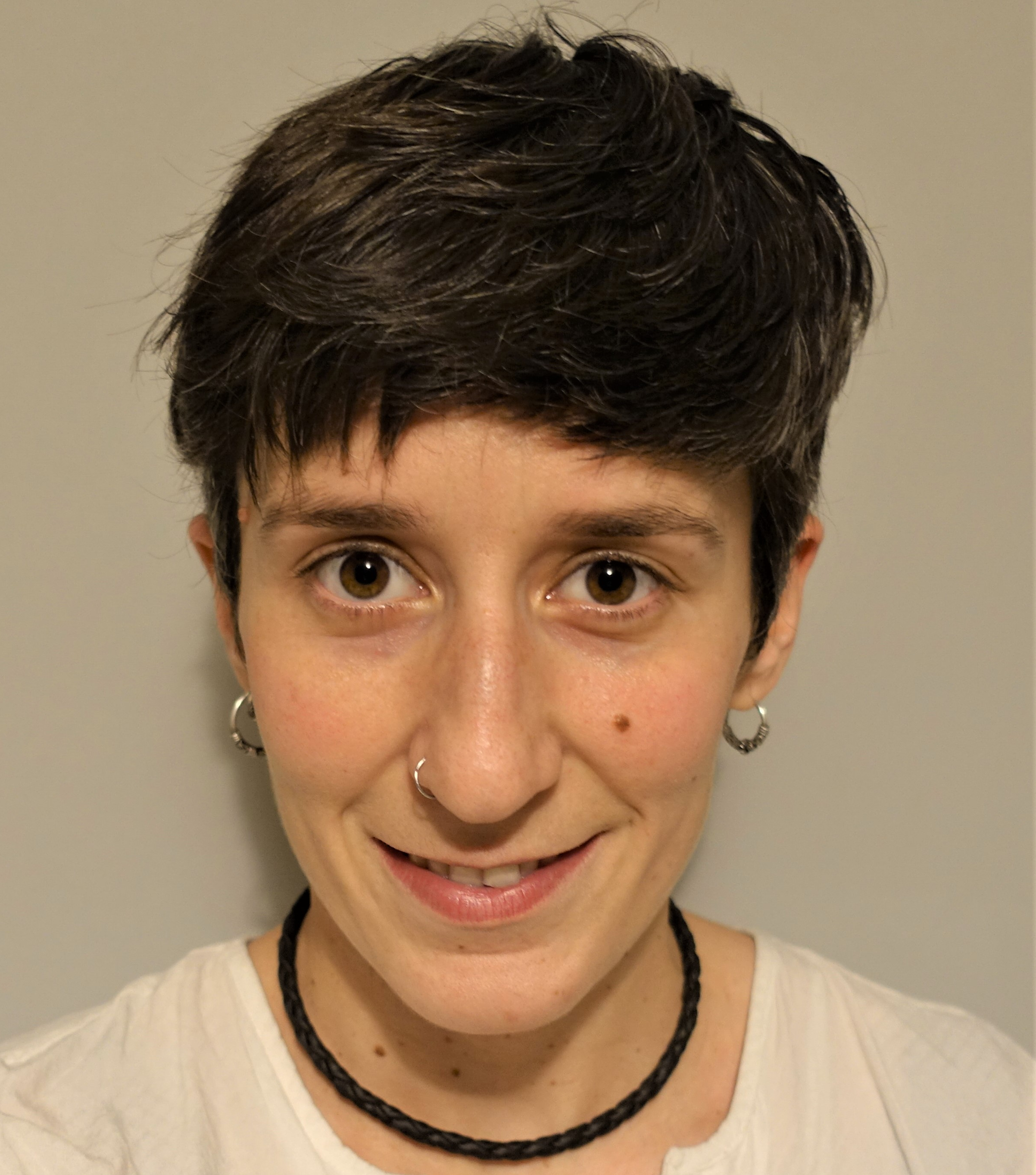}}]{Mar\'{\i}a Cerezo-Maga\~na}
is a Teaching Assistant of Aerospace Engineering at the Universidad Rey Juan Carlos, Madrid, Spain. 
She received the M.Sc. degree in aeronautical engineering from the Universidad Polit\'ecnica de Madrid, Madrid, Spain, and the Ph.D. degree from the Universidad Rey Juan Carlos, in 2015 and 2022, respectively. 
Her research is focused on deterministic and stochastic hybrid optimal control applied to aircraft trajectory optimization.
\end{IEEEbiography}

\begin{IEEEbiography}[{\includegraphics[width=1in,height=1.25in,clip,keepaspectratio]{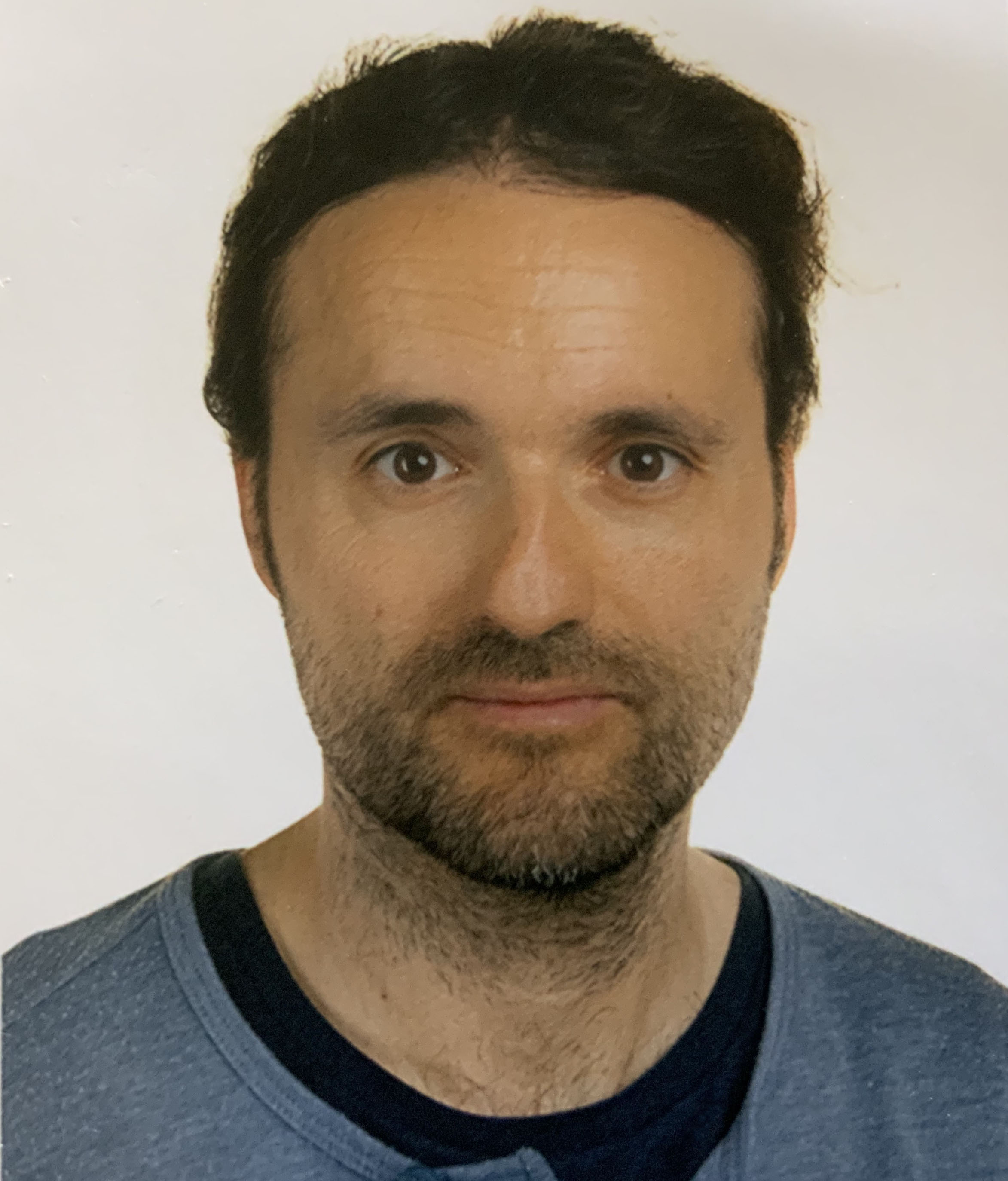}}]{Alberto Olivares}
is a Professor of Statistics and Vector Calculus at the Universidad Rey Juan Carlos, Madrid, Spain. 
He received the M.Sc. degree in mathematics, in 1999, and the B.Sc. degree in statistics, in 2000, both from the Universidad de Salamanca, Salamanca, Spain, and the Ph.D. degree in mathematical engineering from the Universidad Rey Juan Carlos, in 2005.
He was visiting scientist at the Athens University of Economics and Business. His research interests include deterministic and stochastic hybrid optimal control with applications to biomedicine, robotics, aeronautics, and astronautics.
\end{IEEEbiography}

\begin{IEEEbiography}[{\includegraphics[width=1in,height=1.25in,clip,keepaspectratio]{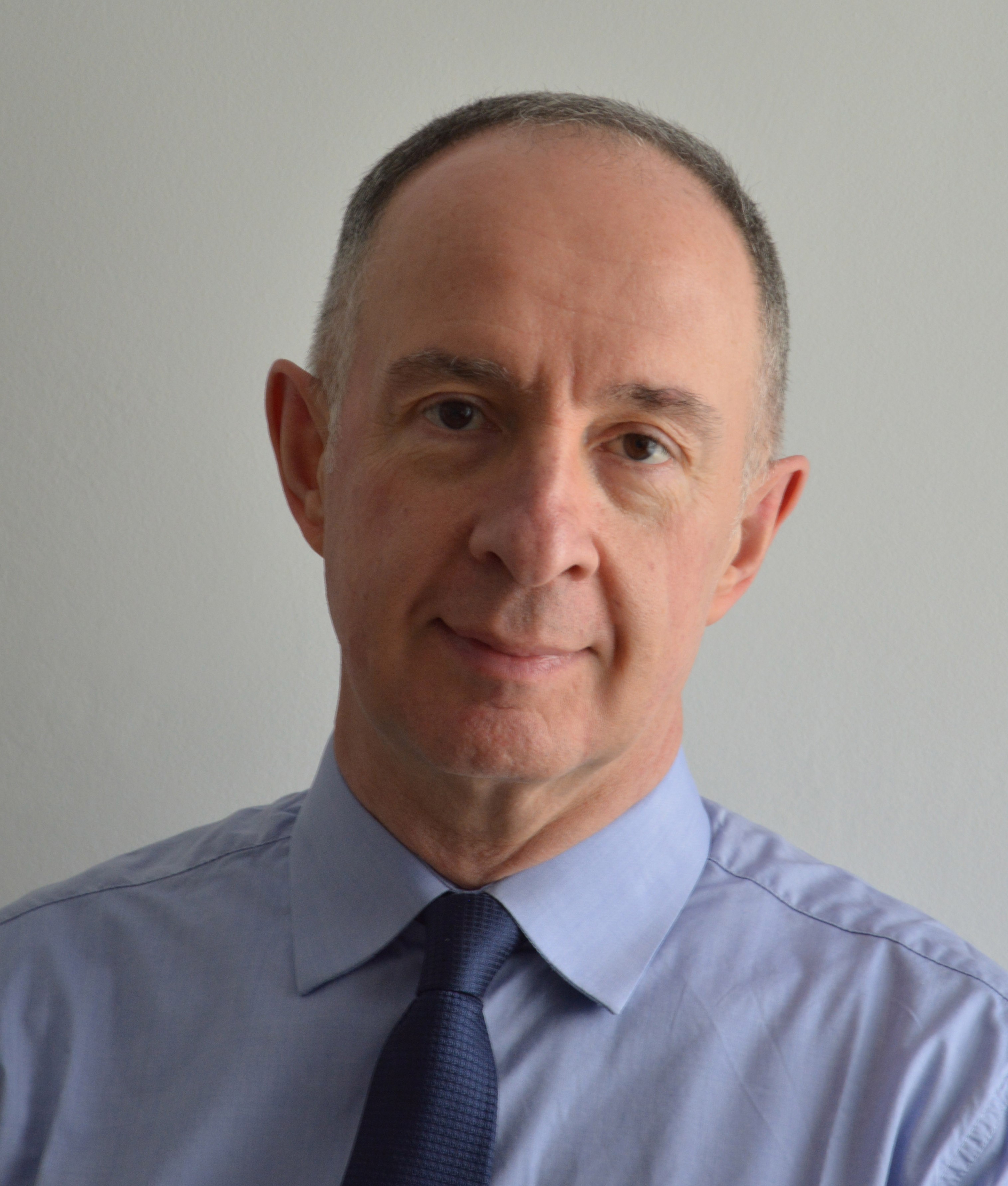}}]{Ernesto Staffetti}
is a Professor of Statistics and Control Systems at the Universidad Rey Juan Carlos, Madrid, Spain. 
He received the M.Sc. degree in automation engineering from the Universit\`a degli Studi di Roma ``La Sapienza," Rome, Italy, and the Ph.D. degree in advanced automation engineering from the Universitat Polit\`ecnica de Catalunya, Barcelona, Spain, in 1995 and 2002, respectively.
He was with the Universitat Polit\`ecnica de Catalunya, Katholieke Universiteit Leuven, Leuven, Belgium, Consejo Superior de Investigaciones Cient\'{\i}ficas, Barcelona, Spain, and University of North Carolina at Charlotte, Charlotte, NC, USA. His research interests include deterministic and stochastic hybrid optimal control with applications to robotics, aeronautics, and astronautics.

\end{IEEEbiography}

\vfill

\end{document}